\documentclass[12pt,reqno]{amsart}
\usepackage{amsmath,amssymb,mathtools,subfigure}
\usepackage{graphicx,color}
\usepackage{tikz,ifthen,cases}
\usetikzlibrary{patterns}
\usepackage{cite}
\usepackage[margin=2.5cm]{geometry}
\usepackage{epsfig}
\usepackage{amsfonts}

\newtheorem{remark}{Remark}[section]

\def\NL{{ M }}  
\def\Nc{{ N }} 
\def\L{{ L }} 
\def\Jump{{ J }} 
\def\Eo{{ \lambda }} 
\def\tt {{ \tau }} 
\def\t {{ t }} 
\def\rhob {\bar{\rho}} 
\def\w{w} 
\def\k{\kappa} 

\def\E{{ \mathbb{E} }} 
\def\rhoc{{\bar{\rho}_{\rm crit}}}
\def\Fbar{{\langle F \rangle}}
\def\vbar{{\langle v \rangle}}

\begin{document}

\title[Nonlocal Traffic Flow Models]{Accelerated Kinetic Monte Carlo methods for general nonlocal traffic flow models}

\author[Y. Sun]{Yi Sun}
\address[Yi Sun]{\newline Department of Mathematics, \
 University of South Carolina, 1523 Greene St., Columbia, SC 29208, USA}
\email{yisun@math.sc.edu}
\author[C. Tan]{Changhui Tan}
\address[Changhui Tan]{\newline Department of Mathematics, \
 University of South Carolina, 1523 Greene St., Columbia, SC 29208, USA}
\email{tan@math.sc.edu}



\date{}
\maketitle

\begin{abstract}
  This paper presents a class of one-dimensional cellular automata (CA) models on traffic flows, featuring nonlocal look-ahead interactions. We develop kinetic Monte Carlo (KMC) algorithms to simulate the dynamics. The standard KMC method can be inefficient for models with global interactions. We design an accelerated KMC method to reduce the computational complexity in the evaluation of the nonlocal transition rates. We investigate several numerical experiments to demonstrate the efficiency of the accelerated algorithm, and obtain the fundamental diagrams of the dynamics under various parameter settings.
\end{abstract}


\section{Introduction}

The mathematical theory on traffic flows has been fast developing in the past century. Many successful models have been proposed, analyzed and simulated
\cite{NER00, CSS00, Hel01, Nag02, Scha02, Ker04, MaD05, KuG11, BeD11, SCN11, TrK13} to understand the emergent phenomena in the traffic networks.
These models can be categorized by different scales.

A famous \emph{macroscopic model} is the Lighthill-Whitham-Richards (LWR)
model \cite{LiW55, Whi74},
\begin{equation}\label{eq:LWR}
  \partial_t\rho+\partial_x(\rho u)=0,\quad u=u_{\max}(1-\rho).
\end{equation}
Here, $\rho$ denotes the normalized density of the traffic, taking values in $[0,1]$. $u$ is the velocity, taking maximum value $u_{\max}$ if $\rho=0$, and becomes 0 if the maximum density $\rho=1$ is reached.
The LWR model can be equivalently expressed as a scalar conservation law:
\begin{equation}\label{eq:LWRflux}
  \partial_t\rho+\partial_x\big(f(\rho)\big)=0, \quad f(\rho)=\rho u(\rho)=u_{\max}\rho(1-\rho),
\end{equation}
where $f$ is the nonlinear flux.
This elegant model captures the wave breakdown phenomenon, which is responsible for the creation of traffic jams.

The LWR model \eqref{eq:LWR} has many extensions. One direction is to consider the \emph{nonlocal slowdown effect}: drivers intend to slow down if heavy traffic is ahead. This would involve a nonlocal look-ahead interaction
\begin{equation}\label{eq:SK}
  \partial_t\rho+\partial_x(\rho u)=0,\quad u=u_{\max}(1-\rho)\exp\left[-\int_0^\infty K(y)\rho(x+y)dy\right],
\end{equation}
where $K$ is the look-ahead kernel. The model was first introduced by Sopasakis and Katsoulakis (SK) in \cite{SKa06}, with
\begin{equation}\label{eq:K-SK}
  K(x)=\begin{cases}1&x\in[0,a],\\0&\text{otherwise.}\end{cases}
\end{equation}
The kernel features a look-ahead distance $a$ and a constant weight.
Another class of kernels was discussed in \cite{CKP14} for pedestrian flows, where
\begin{equation}\label{eq:K-CKP}
  K(x)=\begin{cases}2\big(1-\frac{x}{a}\big)&x\in[0,a],\\0&\text{otherwise,}\end{cases}
\end{equation}
followed by an extensive numerical study.
The wave breakdown phenomenon for the SK model \eqref{eq:SK} and related nonlocal models has been studied in \cite{KuP09, LeL15, Lee20}. Recently, it is shown in \cite{LT22} that the nonlocal slowdown effect can help avoid traffic jams for a family of initial configurations.

Another extension to the LWR model \eqref{eq:LWR} is on the flux in \eqref{eq:LWRflux}. Observe that $f$ is a concave function of $\rho$ with an even symmetry at $\rho=1/2$. This does not agree with the fundamental diagrams from statistical data from real traffic networks, see e.g. \cite{SCN11, KuG11}. A family of fluxes were  introduced in \cite{Pip67} with
\begin{equation}\label{eq:fluxcc}
  f(\rho)=u_{\max}\rho(1-\rho)^\Jump,\quad \Jump>1.
\end{equation}
The fluxes are right-skewed and non-concave, fitting better with the experimental data. The non-concavity can lead to a different type of wave breakdown, as discussed in \cite{Lee19, HT22}.

A general class of traffic flow models with fluxes in \eqref{eq:fluxcc} and nonlocal look-ahead interactions takes the form
\begin{equation}\label{eq:main}
  \partial_t\rho+\partial_x(\rho u)=0,\quad u=u_{\max}(1-\rho)^\Jump g\left(\int_0^\infty K(y)\rho(x+y)dy\right).
\end{equation}
Here $g: [0,+\infty) \to [0,1]$ is the function characterizing the \emph{slowdown factor}. Naturally, $g$ is a nonincreasing function, representing that heavier traffic ahead leads to more slowdown. Also, $g(0)=1$, namely no slowdown if there is no traffic ahead.
A typical choice of $g$ is the Arrhenius relation $g(x)=e^{-x}$. This leads to a generalization of the SK model \eqref{eq:SK}.
Other choices of the function $g$ are $g(x)=(1-x), (1-x)^2$, etc. See for instance the model studied by Bressan and Shen \cite{BrS20}.

\medskip

In this work, we are interested in the \emph{microscopic models} that are closely related to \eqref{eq:main}.

One class of microscopic dynamics is the \emph{agent-based models}, featuring interacting ODE systems on the locations $(x_i)_{i=1}^N$ and/or velocities $(v_i)_{i=1}^N$ of cars. Many models are proposed and studied \cite{GHR61, BHH95, LWS99, Nag99, NSH02}, including the consideration of the look-ahead interactions \cite{WBH04}.

Another class of microscopic dynamics, which is our main concern, is the \emph{lattice models}.
The road is configured as a fixed lattice. Each cell has values $1$ (car is present) or $0$ (car is absent). Explicit rules for car movement on the lattice cells are described to represent the traffic flow. The lattice models, also known as \emph{cellular automata} (CA) models \cite{Wol86, Wol94, Wol02}, have been widely used to represent traffic flows. A vast literature exists addressing various analytical and numerical techniques for models of this type \cite{CrL86, NaS92, BML92, NaP95, BSS98, KSS00, KSS02, LiT05, Nel06}. Compared with the agent-based models, CA models are simpler to implement and are more amenable to numerical investigation.

In \cite{SKa06}, a CA model with Arrhenius type look-ahead interactions was proposed. Through a \emph{semi-discrete mesoscopic stochastic process}, the SK model \eqref{eq:SK} can be formally derived as a coarse-grained hydrodynamic limit of the CA model. See also \cite{HST14} for an improved mesoscopic model that connects the microscopic and macroscopic dynamics. Further extensions include multilane \cite{DuS07}, multiclass \cite{AlS08} and multi-dimensions \cite{SuT14}.

A new class of CA models was introduced by the authors in \cite{SuT20}, recovering the fluxes in \eqref{eq:fluxcc}. A remarkable discovery is that the parameter $\Jump$ represents the number of cells that a car advances in one movement. Following this idea, we describe a large class of CA models in Section \ref{sec:micro}. A formal derivation is provided that connects the CA models to the macroscopic dynamics \eqref{eq:main}, with general choices of the function $g$.

A major focus of this paper is on the numerical implementation of the CA models. One widely used method for vehicular flows and pedestrian flows is the
\emph{Metropolis Monte Carlo (MMC) method} \cite{MRR53}. It is easy to implement, but can be inefficient. Indeed, selected events are sometimes rejected because the acceptance probability is small, in particular when a system approaches the equilibrium, or the car density is high.

To improve computational efficiency, we use the \emph{kinetic Monte Carlo (KMC) algorithm} \cite{BKL75} due to its main feature: \emph{rejection-free}.
Compared to the MMC method, the KMC method requires fewer events to be executed in order to reach a target time, especially when the system is closed to the equilibrium. The KMC method has been successfully applied to traffic flow models \cite{SuT14,SuT20} with special types of nonlocal interactions.

One disadvantage of the KMC method is that the transition rates of all possible events have to be calculated prior to the selection of an event, while the MMC method only requires the transition rate for the selected event. For models with global look-ahead interactions, it is computationally costly to obtain all transition rates due to its nonlocal nature. To overcome such inefficiency, we introduce a new way to evaluate the transition rates,  updating from the previous steps. Taking advantage of the fact that only one car advances in one event, the updates are much cheaper compared with direct evaluations, reducing the cost from $\mathcal{O}(\NL^2)$ to $\mathcal{O}(\NL)$, where $\NL$ denotes the number of cells in the lattice. We call the new procedure the \emph{accelerated KMC method}.

We apply our accelerated KMC method to the nonlocal traffic flow models with a variety of parameter setups. The computational efficiency is verified through numerical experiments. We also obtain the fundamental diagrams of these dynamics, and discuss the relation to the PDE models as coarse-grained limits of the CA models.

The rest of the paper is organized as follows. In Sec.~\ref{sec:micro}, we introduce the CA models with general nonlocal interaction rules, and the connections to macroscopic models like \eqref{eq:main}. In Sec.~\ref{sec:KMC}, we describe the standard KMC algorithm and introduce the new accelerated KMC method. The computational efficiency is also analyzed and compared. In Sec.~\ref{sec:numerics}, we provide a series of numerical simulations to demonstrate the efficiency of our accelerated KMC method. We also generate fundamental diagrams for the CA models and compare them with the macroscopic models. Finally, we state our conclusions in Sec.~\ref{sec:conclusion}.

\section{Cellular Automata Models with Nonlocal Interaction Rules}\label{sec:micro}

In this section, we describe the construction of cellular automata (CA) models for 1D traffic flow, and the connection to the macroscopic models.

The CA models are defined on a periodic lattice $\mathcal{L}$ with $\NL$ evenly spaced cells, $\mathcal{L}=\{1,2,\ldots,\NL\}$. For simplicity,  we assume that all cars move toward one direction on a single-lane loop highway with no entrances or exits. The configuration at each cell $i\in \mathcal{L}$ is defined by an index $\sigma_i$:

\[
\sigma_i=
\begin{cases}
    1 \qquad \mbox{if a car occupies cell $i$}, \\
    0 \qquad \mbox{if the cell $i$ is empty}.
   \end{cases}                                       
\]
The state of the system is represented by $\sigma=\{\sigma_i\}_{i=1}^{\NL}$, which lies in the configuration space $\Sigma=\{0,1\}^\NL$.
We denote $\Nc$ the number of cars. Clearly we have
\[\Nc=\sum_{i=1}^{\NL}\sigma_i.\]

\subsection{Nonlocal interaction rules}\label{subsec:rule}
Car movements can be represented by the transitions in the state of the system, which follow the spin-exchange dynamics \cite{Lig85}: two nearest-neighbor lattice cells exchange values in each transition. Since all cars move to the right, the only possible configuration changes are of the form
\begin{equation}\label{eq:jump1}
\{\sigma_i=1, \sigma_{i+1}=0\} \rightarrow \{\sigma_i=0, \sigma_{i+1}=1\},
\end{equation}
meaning that the car located at the $i$-th cell moves to the $(i+1)$-th cell when it was not occupied. A generalized spin-exchange dynamics introduced in \cite{SuT20} allows the following types of configuration changes
\begin{equation}\label{eq:jump}
\{\sigma_i=1, \sigma_{i+1}=\cdots=\sigma_{i+\Jump}=0\} \rightarrow
\{\sigma_i=\cdots=\sigma_{i+\Jump-1}=0, \sigma_{i+\Jump}=1\}.
\end{equation}
This represents that the car located at the $i$-th cell moves $\Jump$ cells to the right provided that none of $\Jump$ cells in front was occupied. A novel discovery in \cite{SuT20} is that the parameter $\Jump$ determines the macroscopic fluxes \eqref{eq:fluxcc}.

The transition rate for \eqref{eq:jump1} depends on spatial one-sided interactions and a look-ahead feature to represent drivers' behavior. It takes the form
\begin{equation}
    r_i=\frac{\omega_0}{\Jump}g(\w_i).   \label{eq:rate}
\end{equation}
Here, the prefactor $\omega_0=1/\tau_0$ corresponds to the car moving frequency or speed and $\tau_0$ is the characteristic time. The normalization factor $1/\Jump$ makes sure that the estimated velocity is comparable among different choices of $\Jump$. The function $g$ is the same as in the macroscopic dynamics \eqref{eq:main} that describes the slowdown factor.
The quantity $\w_i$ encodes the weighted nonlocal information ahead
\begin{equation}\label{eq:weight}
  \w_i=\frac{1}{\NL}\sum_{j=1}^\NL \k_{j-i}\sigma_j.
\end{equation}
A larger value of $\w_i$ means heavier traffic ahead. The kernel $\{\k_i\}$ is a microscopic analogue of the look-ahead kernel $K$. For instance, for the SK model \eqref{eq:K-SK},
\begin{equation}\label{eq:kiSK}
  \k_i=\begin{cases}1&i=1,\cdots,\L, \\0&\text{otherwise,}\end{cases}
\end{equation}
where $\L=a\NL$ is the microscopic look-ahead distance. In general we shall assume that the kernel is bounded, namely there exists a constant $\overline{K}$ such that \[0\leq \k_i\leq \overline{K}.\] Under the assumption of a looped highway, the kernel $\{\k_i\}$ is $\NL$-periodic. By convention we set $\k_0=0$.

\subsection{Connection to macroscopic models}\label{sec:macro}
In this section, we formally derive the nonlocal traffic models \eqref{eq:main}
from the CA models.
The derivation is a generalization of \cite{SuT20}.

Let us first obtain a semi-discrete mesoscopic model.
In a time step $\Delta\tt$, the probability of the configuration change
\begin{equation}\label{eq:prob}
  \mathbb{P}\big(\{\sigma_i=1, \sigma_{i+1}=\cdots=\sigma_{i+\Jump}=0\} \rightarrow
  \{\sigma_i=\cdots=\sigma_{i+\Jump-1}=0, \sigma_{i+\Jump}=1\}\big)
  =(\Delta\tt)~ r_i,
\end{equation}
where the rate $r_i$ is given in \eqref{eq:rate}.

Define $\sigma(\tt)=\{\sigma_i(\tt)\}_{i=1}^{\NL}$ be
a continuous-in-time stochastic process with a generator
\begin{equation}\label{eq:generator}
(A\psi)(\tt)=\lim_{\Delta\tt\to0}\frac{\E\big[\psi\big(\sigma(\tt+\Delta\tt)\big)\big]
-\psi\big(\sigma(\tt)\big)}{\Delta\tt},
\end{equation}
for any test function $\psi:\Sigma\to\mathbb{R}$, where $\tt$ is the time variable. All possible configuration changes from $\sigma(\tt)$ to $\sigma(\tt+\Delta\tt)$ obey the transition rule \eqref{eq:prob}.
We have
\begin{equation}\label{eq:Egen}
  \frac{d}{d\tt}\E\psi=\E[A\psi].
\end{equation}
In particular, taking $\psi(\sigma)=\sigma_i$,  we can calculate \eqref{eq:generator}
explicitly as follows
\begin{equation}\label{eq:Asigma}
A\sigma_i(\tt)=-r_i(\tt)\sigma_i(\tt)\prod_{j=1}^\Jump\big(1-\sigma_{i+j}(\tt)\big)
+r_{i-\Jump}(\tt)\sigma_{i-J}(\tt)\prod_{j=1}^\Jump \big(1-\sigma_{i-\Jump+j}(\tt)\big)
=: F_{i-\Jump}(\tt)-F_i(\tt),
\end{equation}
where $F_i$ is defined as
\[F_i(\tt)=r_i(\tt)\sigma_i(\tt)\prod_{j=1}^\Jump\big(1-\sigma_{i+j}(\tt)\big).\]

Let $\rho_i(\tt)=\E[\sigma_i(\tt)]=\mathbb{P}\big(\sigma_i(\tt)=1\big)$. Then, from \eqref{eq:Egen} and \eqref{eq:Asigma}, the dynamics of $\{\rho_i\}_{i=1}^{\NL}$ reads
\begin{equation}\label{eq:semirho}
  \frac{d}{d\tt}\rho_i(\tt)=\E[A\sigma_i(\tt)]=\E[F_{i-J}(\tt)]-\E[F_i(\tt)].
\end{equation}
Note that the right hand side of the equation is not yet a closed form of $\{\rho_i(\tt)\}_{i=1}^{\NL}$. We shall approximate the term $\E[F_i(\tt)]$ and make a closure to the system.

We impose the \emph{propagation of chaos} hypothesis, which means that $\{\sigma_i(\tt)\}_{i=1}^{\NL}$ are independent to each other, namely
\begin{equation}\label{eq:chaos}
  \E[\sigma_i(\tt)\sigma_j(\tt)]=\E[\sigma_i(\tt)]~\E[\sigma_j(\tt)],\quad \forall~i\neq j,~~t\geq0.
\end{equation}
Note that due to nonlocal interactions, the hypothesis \eqref{eq:chaos} is not true for a system with fixed $\NL$ cells. However, as the number of cells $\NL$ tends to infinity, the system can become chaotic, and condition \eqref{eq:chaos} can be valid as $\NL\to\infty$.

By formally assuming the chaotic condition \eqref{eq:chaos}, we get
\[\E[F_i(\tt)]= \rho_i(\tt)\prod_{j=1}^J\big(1-\rho_{i+j}(\tt)\big)~ \E\big[~r_i(\tt)~|~\sigma_i(\tt)=1,\sigma_{i+1}(\tt)=\cdots=\sigma_{i+\Jump}(\tt)=0~\big].
\]
For the rest of the section, we drop the $\tt$-dependence for simplicity.

To estimate the rate $r_i$ in \eqref{eq:rate}, we perform a formal Taylor expansion
of $g$ on $\w_i$ around its expectation and get
\begin{equation}\label{eq:Eri}
\E[r_i]=\frac{\omega_0}{\Jump}g\big(\E(\w_i)\big)+\frac{\omega_0}{\Jump}\sum_{n=1}^\infty \frac{g^{(n)}\big(\E(\w_i)\big)}{n!}\,\E\left[\big(\w_i-\E[\w_i]\big)^n\right].
\end{equation}
Note that from \eqref{eq:weight}, $\E[\w_i]$ can be expressed in terms of $\{\rho_i\}_{i=1}^{\NL}$.
\[
  \E[\w_i]=\frac{1}{\NL}\sum_{j=1}^\NL\k_{j-i}\rho_j.
\]
We claim that all terms in the series on the right hand side of \eqref{eq:Eri} vanish as $\NL\to\infty$.
For $n=1$, it is clear that $\E\big[\w_i-\E[\w_i]\big]=0$. For $n=2$, we compute
\begin{align*} &\E\left[\big(\w_i-\E[\w_i]\big)^2\right]=\frac{1}{\NL^2}\E\left[\left(\sum_{j=1}^{\NL}\k_{j-i}(\sigma_j-\rho_j)\right)^2\right]\\& =\frac{1}{\NL^2}\E\left[\sum_{j=1}^{\NL}\k_{j-i}^2(\sigma_j-\rho_j)^2\right]+
\frac{2}{\NL^2}\E\left[\sum_{j=1}^{\NL}\sum_{\ell=1}^{\NL}\k_{j-i}\k_{\ell-i}(\sigma_j-\rho_j)(\sigma_\ell-\rho_\ell)\right].
\end{align*}
By condition \eqref{eq:chaos}, the cross terms
\[\E[(\sigma_j-\rho_j)(\sigma_l-\rho_l)]=\E[\sigma_j-\rho_j]\E[\sigma_l-\rho_l]=0.\]
For the remaining term, since $|\sigma_j-\rho_j|\leq1$, we have $\E\big[(\sigma_j-\rho_j)^2\big]\leq 1$. As $\k_i$ are bounded, we obtain
\[\E\left[\big(\w_i-\E[\w_i]\big)^2\right]=\frac{1}{\NL^2}\sum_{j=1}^{\NL}\k_{j-i}^2\E\left[(\sigma_j-\rho_j)^2\right]\leq \frac{\overline{K}^2}{\NL}\xrightarrow{~\NL\to\infty~}0.\]
Similarly, higher moments vanishes when $\NL\to\infty$.

Plugging back into \eqref{eq:Eri}, we conclude with
\[\E[r_i]=\frac{\omega_0}{\Jump}~g\left(\frac{1}{\NL}\sum_{j=1}^{\NL}\k_{j-i}\rho_j\right)+o(\NL).\]
The conditional expected rate can be represented similarly as
\[\E\big[r_i~|~\sigma_i=1, \sigma_{i+1}=\cdots=\sigma_{i+J}=0\big]=\frac{\omega_0}{\Jump}~g\left(\frac{1}{\NL}\sum_{j\neq i,\cdots,i+\Jump}\k_{j-i}\rho_j\right)+o(\NL).\]
Hence, $\E[F_i]$ can be approximated in terms of $\{\rho_i(\tt)\}_{i=1}^\NL$ as \[
\E[F_i]=\frac{\omega_0}{\Jump}\cdot\rho_i\prod_{j=1}^J\big(1-\rho_{i+j}\big) g\left(\frac{1}{\NL}\sum_{j\neq i,\cdots,i+\Jump}\k_{j-i}\rho_j\right)+o(\NL).
\]

Now, we are ready to derive the coarse-grained PDE model.
Let us rescale the lattice $\mathcal{L}$ into a fixed interval $\mathbb{T}=[0,1]$, where each cell has length $h=1/\NL$.  The $i$-th cell is rescaled to the interval
$[(i-1)h, ih]$.

Define the macroscopic density
$\rho:~\mathbb{T}\times\mathbb{R}_+\to\mathbb{R}$, where
\[\rho(x,\t) = \rho_i(\tt),\quad\text{with }~x=ih,~~\t= \tt h.\]
Letting $h\to0$, we formally obtain the macroscopic flux
\[
  F(x,\t):=J\cdot\lim_{h\to0}\E[F_i(\tt)]=
   \omega_0\rho(x,\t)\big(1-\rho(x,\t)\big)^\Jump
   g\left(\int_{\mathbb{T}} K(y-x)\rho(y,\t)\,dy\right).
 \]
The dynamics of $\rho$ in \eqref{eq:semirho} becomes the following scalar conservation law:
\begin{align*}
  \partial_\t\rho(x,\t)&=\frac{1}{h}\frac{d}{d\tt}\rho_i(\tt)
  =\frac{\E [F_{i-\Jump}(\tt)]-\E [F_i(\tt)]}{h}\\
&\xrightarrow{~h\to0~}\lim_{h\to0}\frac{F(x-\Jump h,\t)-F(x,\t)}{\Jump h}
=-\partial_x(F(x,\t)).
\end{align*}

We end up with the following coarse-grained PDE model:
\begin{equation}    \label{eq:PDE}
    \partial_t\rho+\partial_x\left(\omega_0\rho(1-\rho)^\Jump
      g\left(\int_{\mathbb{T}}K(y-x)\rho(y,t)dy\right)\right)=0.
\end{equation}
It is the periodic version of the macroscopic model \eqref{eq:main} with
$u_{\max}=\omega_0$.


\section{The Accelerated Kinetic Monte Carlo Method} \label{sec:KMC}

In this section, we focus on the numerical implementation of the CA models in Section \ref{sec:micro}. We use the kinetic Monte Carlo (KMC) method \cite{BKL75} to simulate the spin exchange dynamics. Compared with the Metropolis Monte Carlo (MMC) method \cite{MRR53}, the KMC method has a major advantage: \emph{rejection-free}. In each step, the transition rates for all possible changes from the current configuration are calculated and then a new configuration is chosen with a probability proportional to the rate of the corresponding transition. The other feature of the KMC method is its capability of providing a more accurate description of the real-time evolution of a traffic system in terms of these transition rates since the KMC method is more suitable for simulating the non-equilibrium system.

\subsection{The KMC algorithm}
Let us describe the KMC method that we use for the CA models with nonlocal look-ahead interactions.

We start with some notations. Let us denote  $(i_j)_{j=1}^{\Nc}$ the ordered locations of the cells that are occupied by cars. Here we recall $\Nc$ is the number of cars. Note that all events happen at these occupied cells. According to \eqref{eq:jump}, the event $k$ is that the $k$-th car located at the cell $i_k$ moves $\Jump$ cells to the right, with the rate $r_{i_k}$. The KMC algorithm is built on the assumption that the model features $\Nc$ independent Poisson processes (corresponding to $\Nc$ moving cars on the lattice) with transition rates $r_{i_j}$ in \eqref{eq:rate} that sum up to give the total rate $R=\sum^{\Nc}_{j=1}r_{i_j}$. In each round of the KMC algorithm, we need to do the following.

{\it The KMC algorithm:}

Step 1:  Generate a random number $\xi_1$ from the uniform distribution in $[0, 1]$. Decide which event will take place by using a binary search to choose the event
$k$ such that
\begin{equation} \label{eq.bk1}
  \sum^{k-1}_{j=1} \frac{r_{i_j}}{R} < \xi_1 \leq \sum^k_{j=1} \frac{r_{i_j}}{R},\quad R=\sum_{j=1}^\Nc r_{i_j}.
\end{equation}

Step 2: Check if there are enough vacant cells ahead of the $k$-th car of the selected event. If ``Yes", perform Step 3. If ``No", skip Step 3 and advance directly to Step 4.

Step 3: Perform the selected event (the $k$-th car move $\Jump$ cells to the right) leading to a new configuration. Update the location of the $k$-th car $i_k^{\text{new}}=i_k^{\text{old}}+\Jump$. Also update the total rate $R$ and any rate $r_{i_j}$ that may have changed due to this move.

Step 4: Use $R$ and another random number $\xi_2 \in (0, 1)$ to decide the time it takes for that event to occur (the transition time), i.e., the nonuniform time step $\Delta t=-\ln(\xi_2)/R$, from the exponential distribution described by the rate $R$.

\hfill $\Box$

\begin{remark}
When the system evolves from the initial state to the equilibrium state, the total rate $R$ decreases. From Step 4, we can see that the average of the time step $\Delta t$ will increase with the decreasing $R$ so that the KMC simulation will reach the preset final time sooner. This is another computational advantage over the MMC, which usually sets the time step $\Delta t$ to be a small constant.
\end{remark}

\subsection{List-based KMC methods}
In simulations with a finite number of distinct processes, it is more efficient to consider the groups of events according to their rates \cite{BBS95, Sch02, SCE11}. This is known as the \emph{list-based methods}.

In the context of traffic flow models, the list-based KMC algorithm has been successfully implemented in \cite{SuT14,SuT20} for the SK model where the kernel takes a special form \eqref{eq:kiSK} with a constant value in the look-ahead distance $\L$. In this case, the rate $r_i$ in \eqref{eq:rate} can only take $(\L+1)$ different values, since $w_i$ in \eqref{eq:weight} can only be $0, \frac{1}{\NL},\ldots, \frac{\L}{\NL}$. To speed up the event search process in Step 1 of the KMC algorithm, $(\L+1)$ lists are created. Each list is a collection of all events with the same rate. To find an event, we first perform a binary search in the level of lists, namely searching for list $l$ such that
\begin{equation} \label{eq.bk2}
  \sum^{l-1}_{j=1} \frac{n_jr_j}{R} < \xi_1 \leq \sum^l_{j=1} \frac{n_jr_j}{R},\quad R=\sum_{j=1}^{\L+1} n_jr_j.
\end{equation}
Here, we denote $(r_j)_{j=1}^{\L+1}$ the different rates, and $n_j$ the number of events in the list with rate $r_j$, which is called the {\it multiplicity}. Once a list is chosen, we then randomly pick an event in the list to proceed. The binary search in \eqref{eq.bk2} would cost $\mathcal{O}(\log_2\L)$ operations, smaller than $\mathcal{O}(\log_2\Nc)$ that is needed in \eqref{eq.bk1} when $\L$ is much smaller than $\Nc$.

However, the list-based KMC method does not have an advantage when the rates take many different values. For instance, if the kernel is a discrete version of \eqref{eq:K-CKP}, namely
\begin{equation}\label{eq:kiCKP}
  \k_i=\begin{cases}2\big(1-\frac{i-1/2}{\L}\big)&i=1,\cdots,\L,\\0&\text{otherwise},\end{cases}
\end{equation}
the rate can take $\mathcal{O}(\L^2)$ different values. For general kernel $(\k_i)$, the rates can take as many as $\NL^2$ different values. Since there are a total of $\Nc$ events, many lists would have only one or no event. This makes the list-based method inefficient. Therefore, in the present work, we do not use the list-based KMC method.

\subsection{An accelerated KMC method}
Another major issue is the computational cost of updating the rates $r_{i_j}$ in Step 3 of the KMC algorithm.

For the SK model with a constant kernel \eqref{eq:kiSK} in a look-ahead distance $\L$, the configuration change \eqref{eq:jump} at cell $i_k$ only alter the values of $\sigma_{i_k}$ and $\sigma_{i_k+\Jump}$. It is easy to check that the change may only affect the rates of at most $\Jump$ cars in certain cells from the location $i_k-\L$ to $i_k-\L+\Jump-1$ (if there are cars there). Together with the $k$-th car, we only need to update at most $\Jump+1$ rates. Calculating the new rates from \eqref{eq:rate}-\eqref{eq:weight} costs $\mathcal{O}(\L)$ operations. Since $\Jump=\mathcal{O}(1)$, the total cost for Step 3 is $\mathcal{O}(\L)$.

However, for a general kernel $(\k_i)$ with a look-ahead distance $\L$, e.g. \eqref{eq:kiCKP}, the configuration change \eqref{eq:jump} at cell $i_k$ can affect the rates of all the cars behind the moving $k$-th car in the range of $\L$, from the location $i_k-\L$ to $i_k-1$ (if there are cars there). Together with the $k$-th car, we need to update up to $L+1$ rates. Computing these rates from \eqref{eq:rate}-\eqref{eq:weight} would cost up to a total of $\mathcal{O}(\L^2)$ operations. For a global kernel $(\k_i)$ when $\L=\NL$, all rates $(r_{i_j})_{j=1}^{\Nc}$ can be affected, and the total computational cost can be $\mathcal{O}(\Nc\NL)$.

To accelerate the KMC algorithm, we introduce a new method, aiming to reduce the cost of computing the rates of the $(\Nc-1)$ not-moving cars $(r_{i_j})_{j=1, j\neq k}^{\Nc}$ from $\mathcal{O}(\Nc\NL)$ down to $\mathcal{O}(\Nc)$. The main idea is to calculate the new weights $w_{i_j}$ not from \eqref{eq:weight}, but from the old weights in the previous step prior to the event. More precisely, let us denote $\{\sigma_i^{\text{old}}\}_{i=1}^{\NL}$ and $\{\sigma_i^{\text{new}}\}_{i=1}^{\NL}$ the configurations of the system before and after an event \eqref{eq:jump} located at cell $i_k$, respectively. Clearly, we have
\[\sigma_j^{\text{new}}=\sigma_j^{\text{old}},\quad\forall~j\neq i_k, i_k+\Jump.\]
The only differences are
\[\sigma_{i_k}^{\text{old}}=1, \sigma_{i_k}^{\text{new}}=0,\quad\text{and}\quad\sigma_{i_k+\Jump}^{\text{old}}=0, \sigma_{i_k+\Jump}^{\text{new}}=1.\]
We apply \eqref{eq:weight} and obtain the relation
\begin{align}
  \w_{i_j}^{\text{new}}&=\frac{1}{\NL}\sum_{\ell=1}^{\NL}\k_{\ell-i_j}\sigma_\ell^{\text{new}}
=\frac{1}{\NL}\left(\sum_{\ell\neq i_k, i_k+\Jump}^{\NL}\k_{\ell-i_j}\sigma_\ell^{\text{new}}+\k_{i_k+\Jump-i_j}\right)\label{eq:update}\\
  &=\frac{1}{\NL}\left(\sum_{\ell\neq i_k, i_k +\Jump}^{\NL} \k_{\ell-i_j}\sigma_\ell^{\text{old}}+ \k_{i_k-i_j}\right) +\frac{\k_{i_k+\Jump-i_j}-\k_{i_k-i_j}}{\NL} =\w_{i_j}^{\text{old}}+\frac{\k_{i_k+\Jump-i_j}-\k_{i_k-i_j}}{\NL},\nonumber
\end{align}
for any $j$-th car with $j=1,\ldots,\Nc$ and $j\neq k$. The relation \eqref{eq:update} allows us to obtain $\w_{i_j}^{\text{new}}$ from $\w_{i_j}^{\text{old}}$ using $\mathcal{O}(1)$ operations. Then we can compute the rates $(r_{i_j}^{\text{new}})_{j=1, j\neq k}^{\Nc}$ from \eqref{eq:rate} in $\mathcal{O}(\Nc)$ operations. Note that the relation \eqref{eq:update} is not applicable to the $k$-th car, as its location changes during the event. We shall still update $w_{i_k}$ from \eqref{eq:weight}.  Overall, the total updating cost in Step 3 is reduced to $\mathcal{O}(\Nc+\L)\leq\mathcal{O}(\NL)$. It is a big improvement compared with $\mathcal{O}(\Nc\NL)$.

We shall comment that at the beginning of a KMC simulation, we need to compute the weights $(\w_{i_j})_{j=1}^{\Nc}$ of all cars from \eqref{eq:weight} for initialization, and create a data segment to store these weights so that we can update them by using \eqref{eq:update} on the fly in the subsequent steps.


\section{Numerical experiments}\label{sec:numerics}

In this section, we investigate 1D nonlocal traffic flow models in various parameter regimes with the accelarated KMC method presented in the previous section.

Following \cite{SKa06, SuT14, SuT20}, we set the actual physical length of each cell to $22$ feet ($\approx 6.7$m), which allows for the average car length plus safe distance. Therefore, $1$ mile ($=5280$ feet $\approx 1609$m) is equivalent to $240$ cells. For a car which has average speed of $60$ miles per hour ($\approx 26.8$ m/s), an estimate of time to cross a cell is given by
\[
\Delta \tt_{\rm cell}=\frac{22 \ {\rm feet}}{60 \ {\rm miles/h}} = \frac{1 \ {\rm cell} \times 3600 \ {\rm s}}{60 \times 240 \ {\rm cells}} = \frac{1}{4} {\rm s}.
\]
Therefore, in the KMC simulations we set the characteristic time $\tau_0=0.25$s, and then $\omega_0=4$s$^{-1}$. We mention that other values of $\tau_0$ may be chosen to adjust our model for considering different standards in other regions or countries.

One important group of statistical features that characterize the CA models is the \emph{fundamental diagrams}. Define the average flow $\Fbar$ to be the number of cars passing a fixed detector site per unit time \cite{May90}. This quantity can be measured in a real traffic system. In our study, we run KMC simulations for different initial averaged car densities $\rhob$. Each simulation is run for sufficient long time so the dynamics reaches a stable equilibrium: macroscopically speaking $\rho(x,t)\equiv\rhob$.
$\Fbar$ is taken as a long-time average for each simulation, and is further averaged among several simulations with the same $\rhob$.
The function
\[\Fbar = \Fbar(\rhob)\]
is known as the density-flow fundamental diagram of the corresponding CA model. It is closely related to the flux $f=f(\rho)$ for the coarse-grained PDE models, e.g. \eqref{eq:LWRflux} and \eqref{eq:fluxcc}, when there is no nonlocal interaction. Another statistical quantity is the ensemble-averaged velocity $\vbar$, representing the velocity averaged among all cars and in a long time period.  One can generate density-velocity and flow-velocity diagrams in a similar way.

\subsection{The KMC acceleration}
Our first example demonstrates the computational efficiency of the accelerated KMC method, in comparison with the standard KMC method.

\begin{figure}[!ht]
\begin{center}
\includegraphics[width=.48\textwidth]{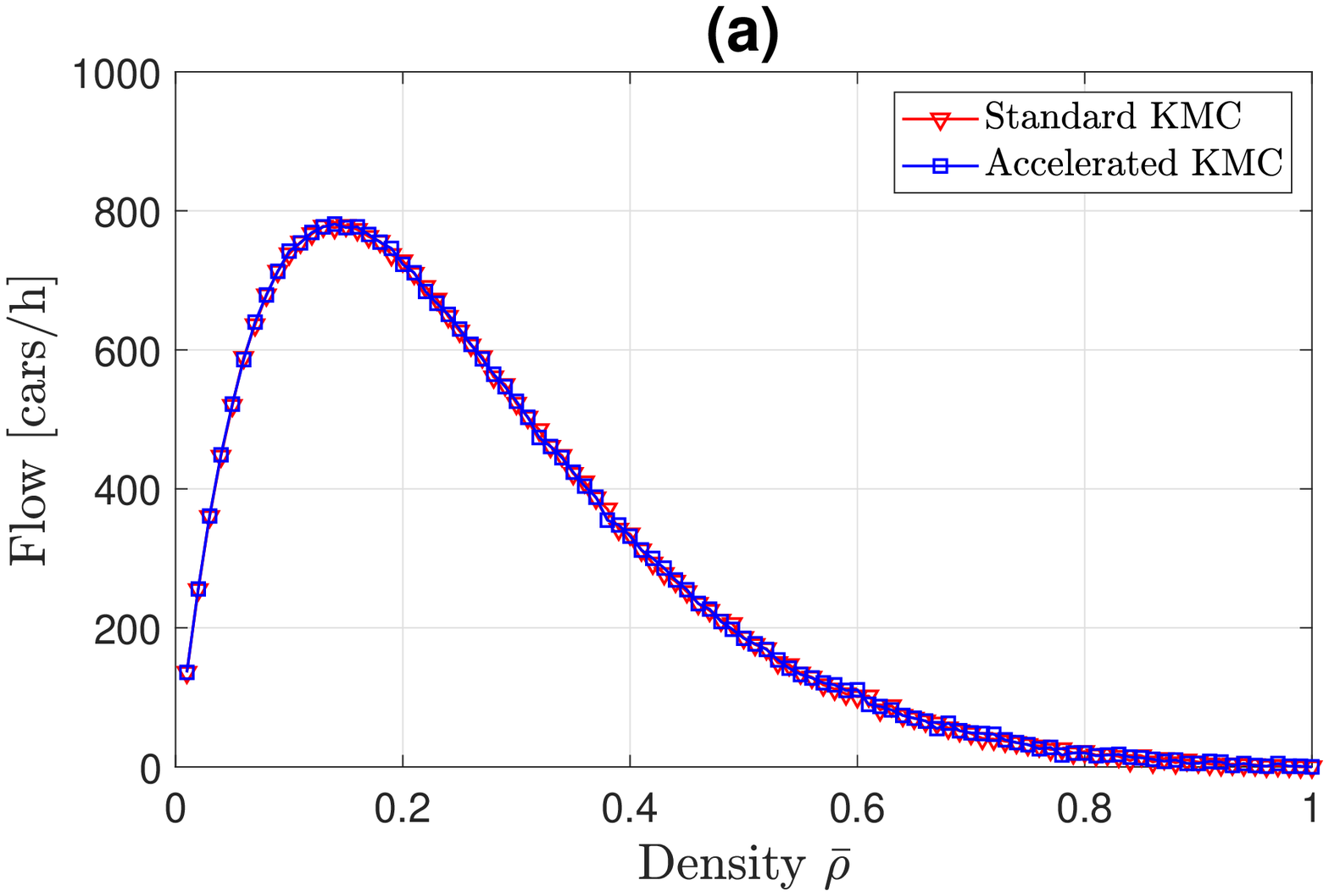} \hfill
\includegraphics[width=.46\textwidth]{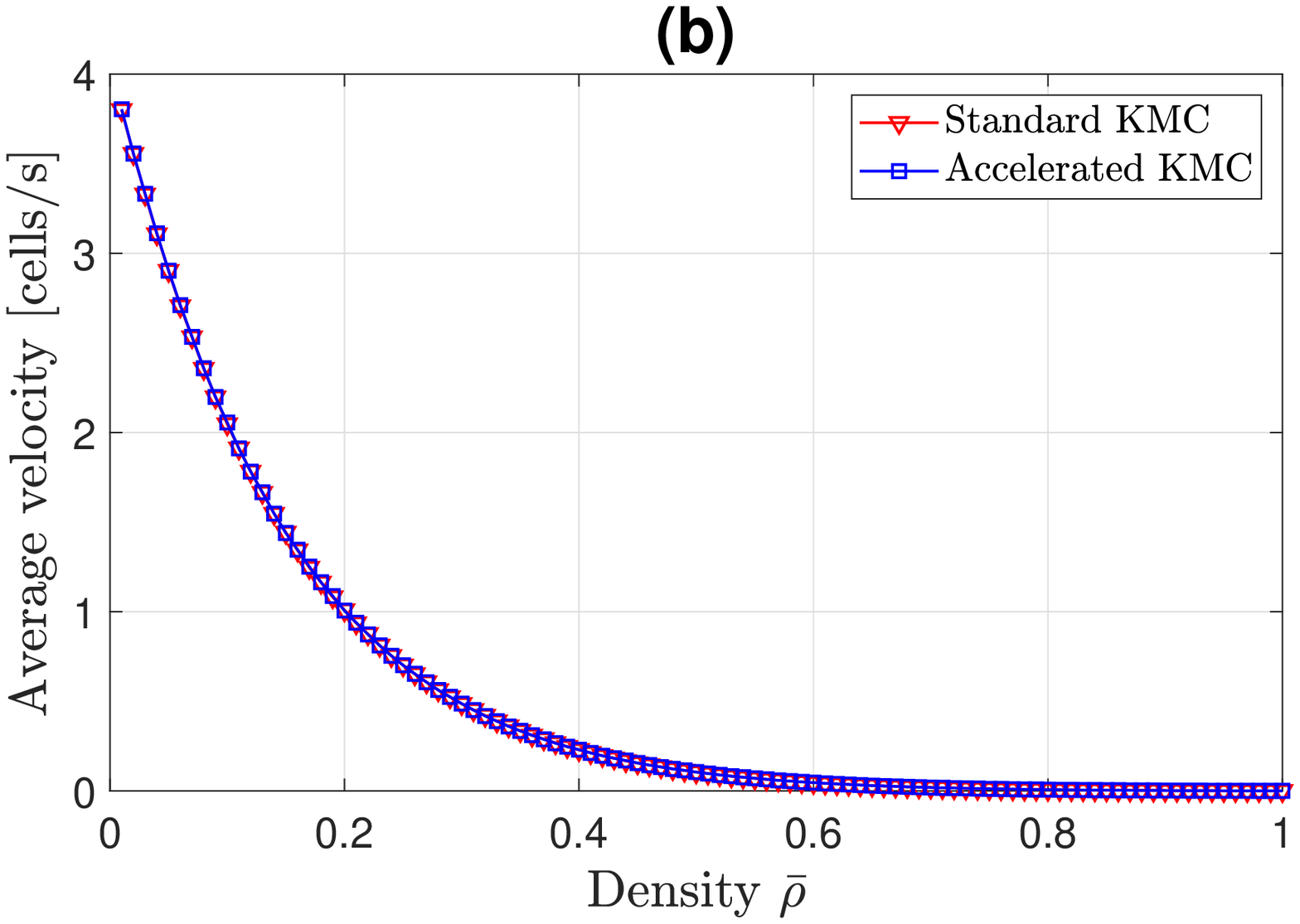}

\vspace{0.2cm}

\includegraphics[width=.47\textwidth]{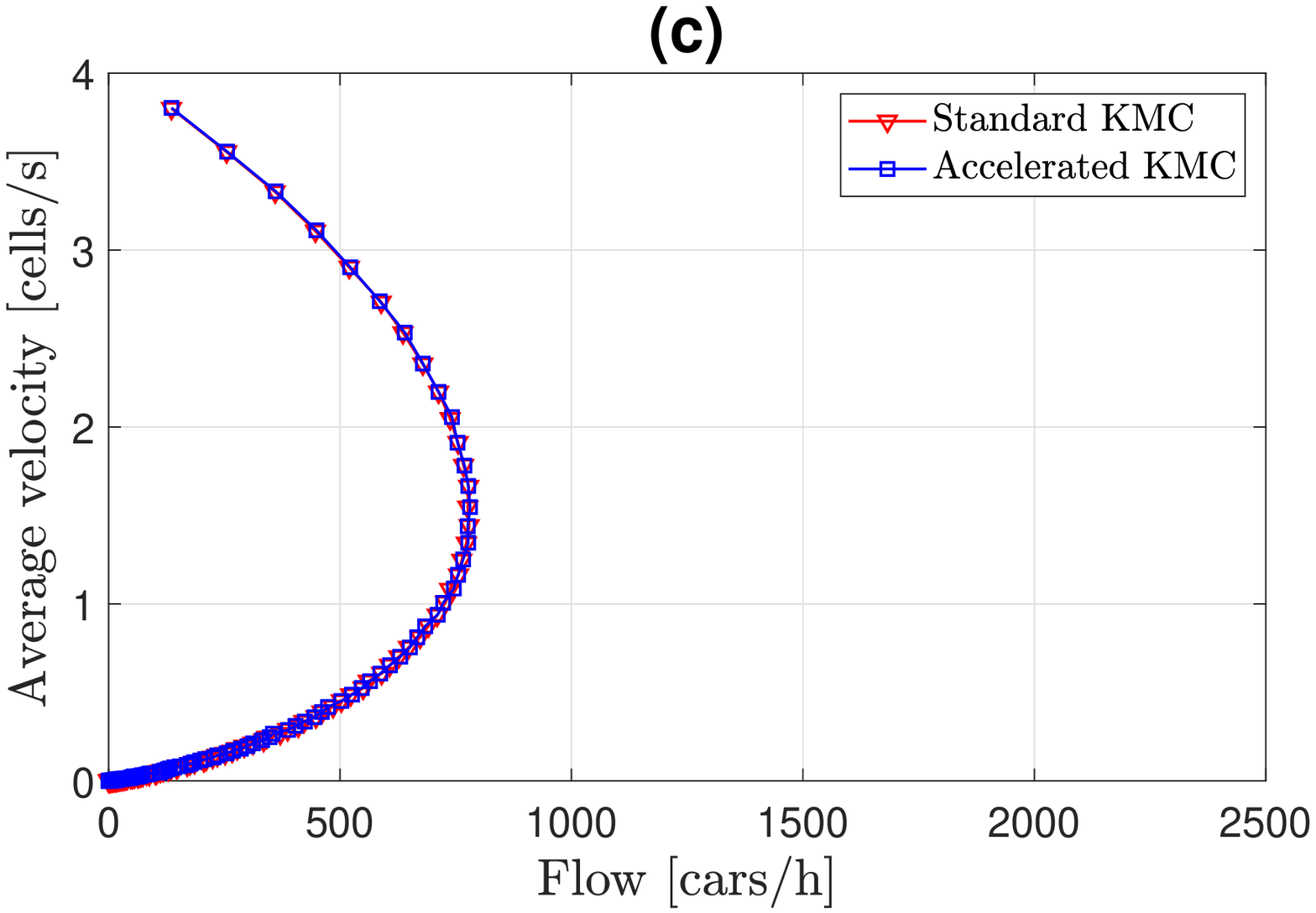}

\caption{Comparison results between the standard KMC and the accelerated KMC algorithms. We take the highway distance of $\approx 2.09$ miles ($\approx 3352$ m, $\NL=500$ cells),  the look-ahead parameter of $\L=500$, the multiple move parameter of $\Jump=1$, and the final time of $1$ h.  (a): Long-time averages of the density-flow relationship; (b): Ensemble-averaged velocity of cars versus the density $\bar{\rho}$; (c): Long-time averages of the flow-velocity relationship.}
                                                           \label{fig:compare}
\end{center}
\end{figure}

Let us consider a loop highway of $\approx 2.09$ miles ($\approx 3352$m, $\NL=500$ cells). Take an Arrhenius type interaction relation $g(x)=e^{-3x}$, and a linearly decay kernel \eqref{eq:kiCKP} with a large look-ahead distance $\L=\NL=500$ for a global interaction. The multiple move parameter is set to be $\Jump=1$. We establish the fundamental diagrams of the density-flow, density-velocity and flow-velocity relationships using the standard KMC and accelerated KMC methods. In particular, we generate random initial distributions with the averaged car density $\rhob$ increasing incrementally from $\rhob=0.01$ to $\rhob=0.99$.
For each $\rhob$, we run 10 simulations with different random number seeds for a long time (1 hour) to get $\Fbar$ and $\vbar$. Fig.~\ref{fig:compare} shows the fundamental diagrams generated from the standard KMC method and the accelerated KMC method. The results agree with each other very well.

\begin{figure}[!ht]
\begin{center}

\includegraphics[width=.48\textwidth]{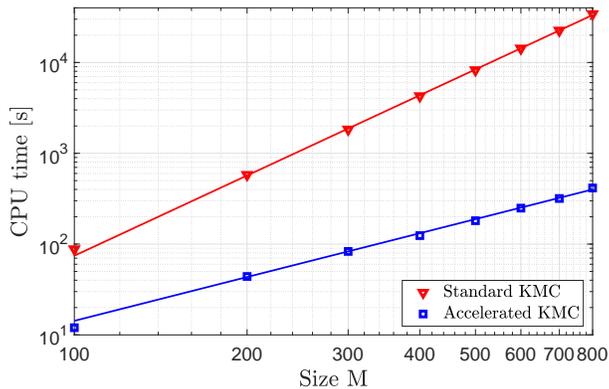}

\caption{A log-log plot of the CPU times with the two KMC algorithms for different highway distance $\NL=100$ to $800$ cells;  power-law lines with exponent $\approx 2.94$ (red) and exponent $\approx 1.59$ (blue) provided for comparison.}
                                                           \label{fig:compareCPU}
\end{center}
\end{figure}

To compare the computational efficiencies of the two algorithms, we now vary the size of the lattice, taking $\NL=100$ to $800$ cells.
Table~\ref{tab:compare} shows the CPU times of computing the simulations described above using both methods.
Fig.~\ref{fig:compareCPU} displays a power-law relationship between the CPU time and the highway distance $\NL$ with a power-law exponent $\approx 2.94$ for the standard KMC algorithm, and a much smaller exponent $\approx 1.59$ for the accelerated KMC algorithm. One can clearly observe the acceleration from our new KMC algorithm.

\begin{center}
\begin{table}
\begin{tabular}{|| r | r | r | r | r | r | r | r | r ||}
\hline
$\NL$ (cells) & $100$ & $200$  & $300$ & $400$ & $500$ & $600$ & $700$ & $800$ \\
\hline\hline
Standard KMC & 93 s & 580 s & 1889 s & 4270 s & 8262 s & 14302 s & 22728 s & 34353 s\\
\hline
Accelerated KMC & 12 s &  44 s &   83 s &   124 s &   181 s &   249 s &  318 s &  415 s \\
\hline \hline
\end{tabular}
\vspace{0.2cm}
\caption{Comparison of CPU times between the standard KMC and the accelerated KMC algorithms for different highway distance $\NL=100$ to $800$ cells, which are also displayed by a log-log plot in Fig.~\ref{fig:compareCPU}.}
\label{tab:compare}
\end{table}
\end{center}

\subsection{A family of models with global interactions}
\label{subsect.strength}

Now we apply the accelerated KMC method to a class of 1D nonlocal traffic flow models with global look-ahead interaction kernels.
The coarse-grained PDE dynamics has the form \eqref{eq:PDE} with
\begin{equation}\label{eq:g-1}
g(x) = \begin{cases}1-x& x\in[0,1],\\0&\text{otherwise},\end{cases}
\end{equation}
and the kernel
\begin{equation}\label{eq:K-BS}
  K^\Eo(x) = \begin{cases}\Eo e^{-\Eo x}&x\geq0,\\0&x<0,\end{cases}
\end{equation}
which is parameterized by $\Eo$.
This family of PDE models was proposed and analytically studied in \cite{BrS20}.

For the CA model, we take the loop highway of $\approx4.17$ miles ($\approx 6704$m, $\NL=1000$ cells).
The discrete analogue of \eqref{eq:K-BS} in the periodic domain that we will use takes the form
\begin{equation}\label{eq:kiBS}
  \k_i^\Eo = \frac{\NL(e^{\Eo/\NL}-1)}{1-e^{-\Eo}}\cdot e^{-\Eo(i/\NL)},\quad i=1,\ldots,\NL.
\end{equation}

Let us focus on the model under two extreme choices of $\Eo$.
First, when $\Eo$ is close to 0, we get from \eqref{eq:kiBS} that
\[\lim_{\Eo\to0}\k_i^\Eo=1,\quad i=1,\ldots,\NL.\]
The kernel becomes uniform. Formally, the macroscopic flux becomes
\[f(\rho)=\omega_0\rho(1-\rho)^\Jump g(\rhob).\]
Therefore, as the dynamics reaches the equilibrium state $\rho(x)\equiv\rhob$, we should expect that the long-time averaged flux and the ensemble-average velocity satisfy
\begin{equation}\label{eq:fluxavg0}
  \Fbar(\rhob)=\omega_0\rhob(1-\rhob)^{\Jump+1},
  \qquad \vbar(\rhob)=\omega_0(1-\rhob)^{\Jump+1}.
\end{equation}
Second, when $\Eo$ approaches $+\infty$, we get from \eqref{eq:kiBS} that
\[\lim_{\Eo\to\infty}\k_i^\Eo=\begin{cases}\NL,& i=1,\\0,&\text{otherwise},\end{cases}\]
Apply this to \eqref{eq:weight} and we have
\[w_i^\Eo=\frac{1}{\NL}\sum_{j=1}^\NL \k_{j-i}\sigma_j=\sigma_{i+1}=0.\]
The last equality holds for any event in \eqref{eq:jump}. Therefore, there is no slowdown ($g(0)=1$), and hence the averaged flux and the averaged velocity satisfy
\begin{equation}\label{eq:fluxavginf}
  \Fbar(\rhob)=\omega_0\rhob(1-\rhob)^{\Jump},
  \qquad \vbar(\rhob)=\omega_0(1-\rhob)^{\Jump}.
\end{equation}
We would like to point out that for the macroscopic dynamics \eqref{eq:PDE} with \eqref{eq:K-BS}, $K^\Eo$ converges to a Dirac delta as $\Eo\to\infty$, and the limiting PDE reads
\[\partial_t\rho+\partial_x(\omega_0\rho(1-\rho)^\Jump g(\rho))=0.\]
The macroscopic flux of the limiting system does not match with \eqref{eq:fluxavginf}. This reveals an interesting effect: the two limits $\NL\to\infty$ and $\Eo\to\infty$ do not commute.

\begin{figure}[p]
\begin{center}
\includegraphics[width=.48\textwidth]{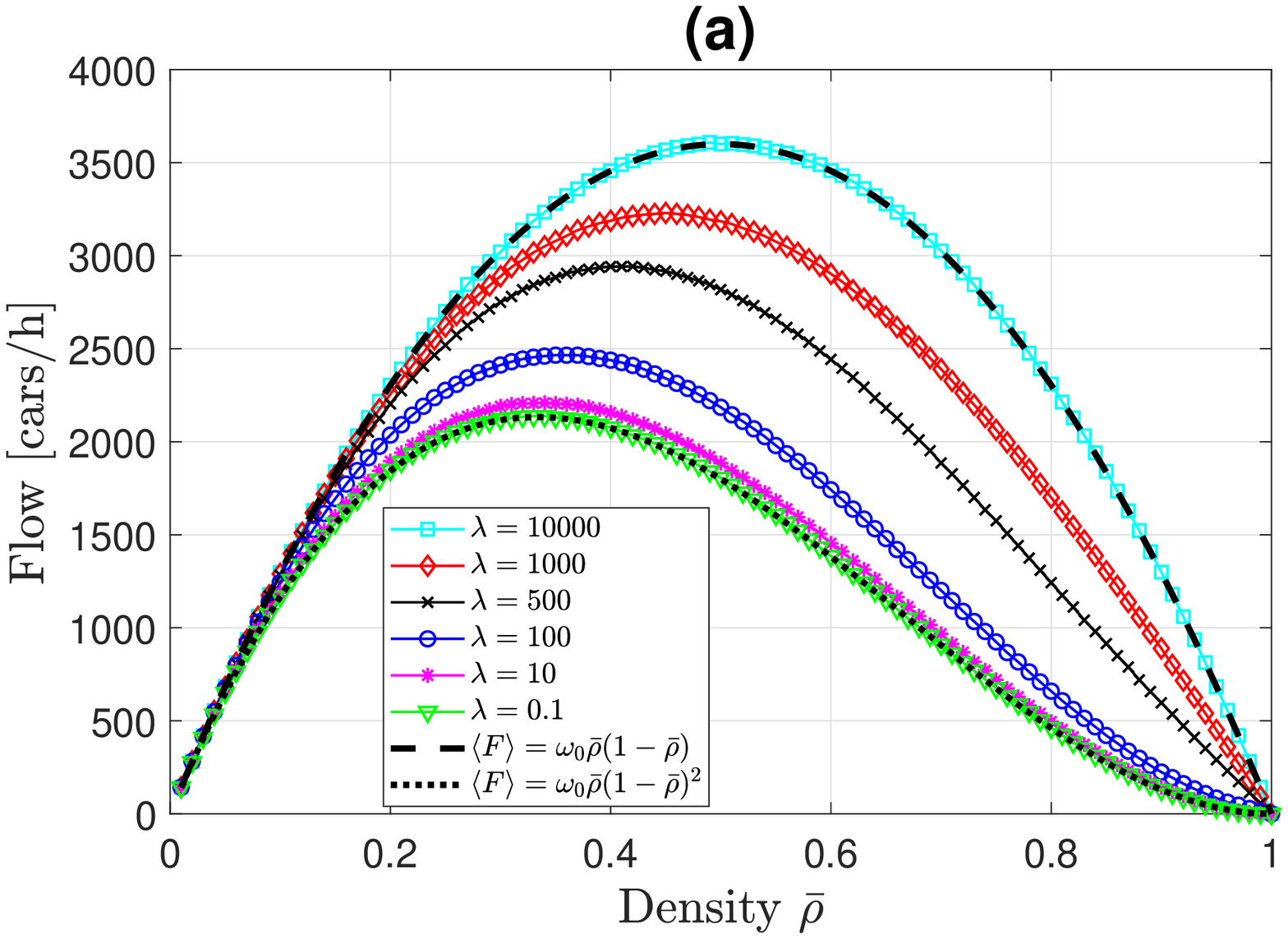} \hfill
\includegraphics[width=.48\textwidth]{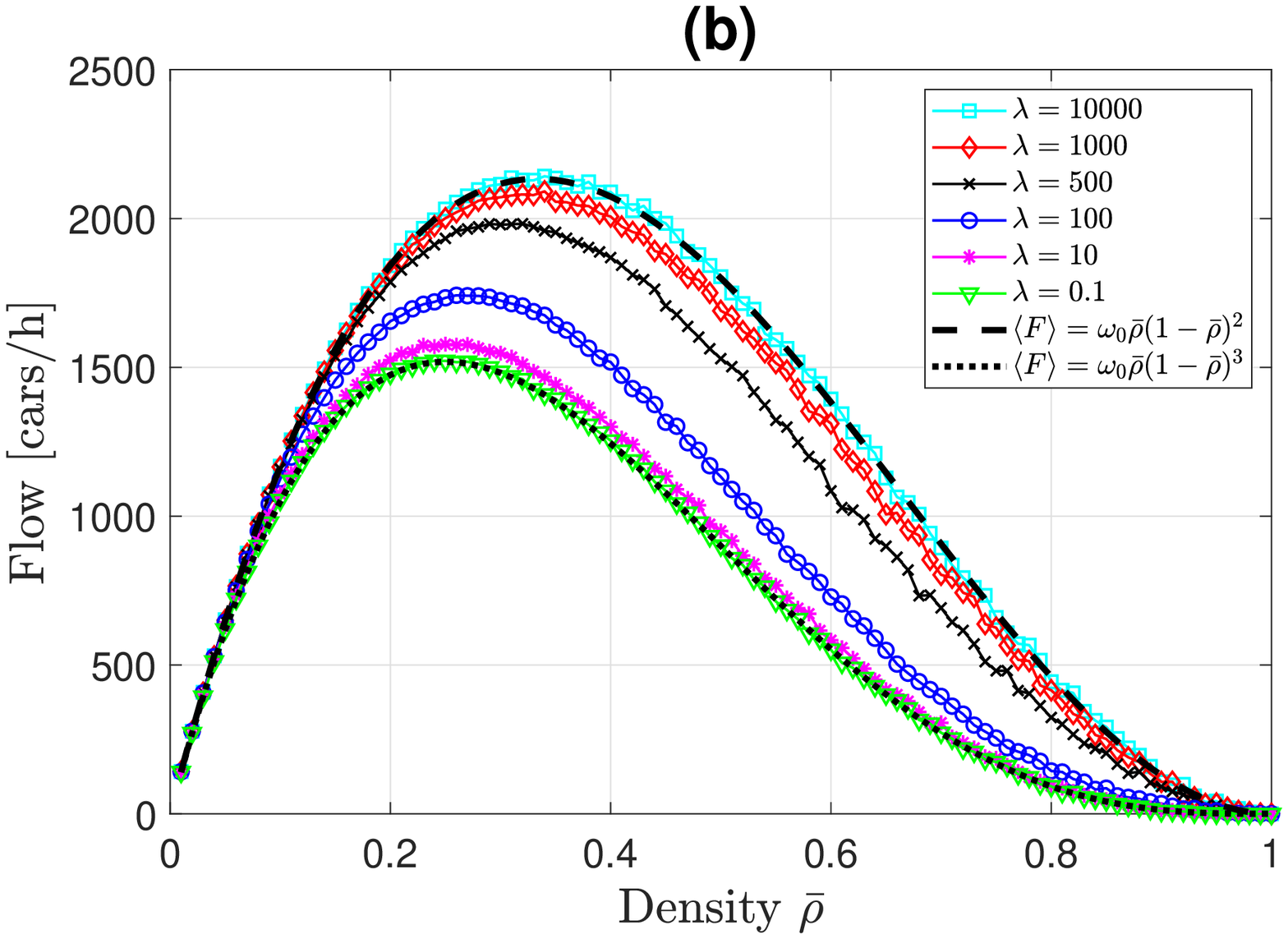}

\vspace{0.2cm}

\hspace{0.2cm} \includegraphics[width=.46\textwidth]{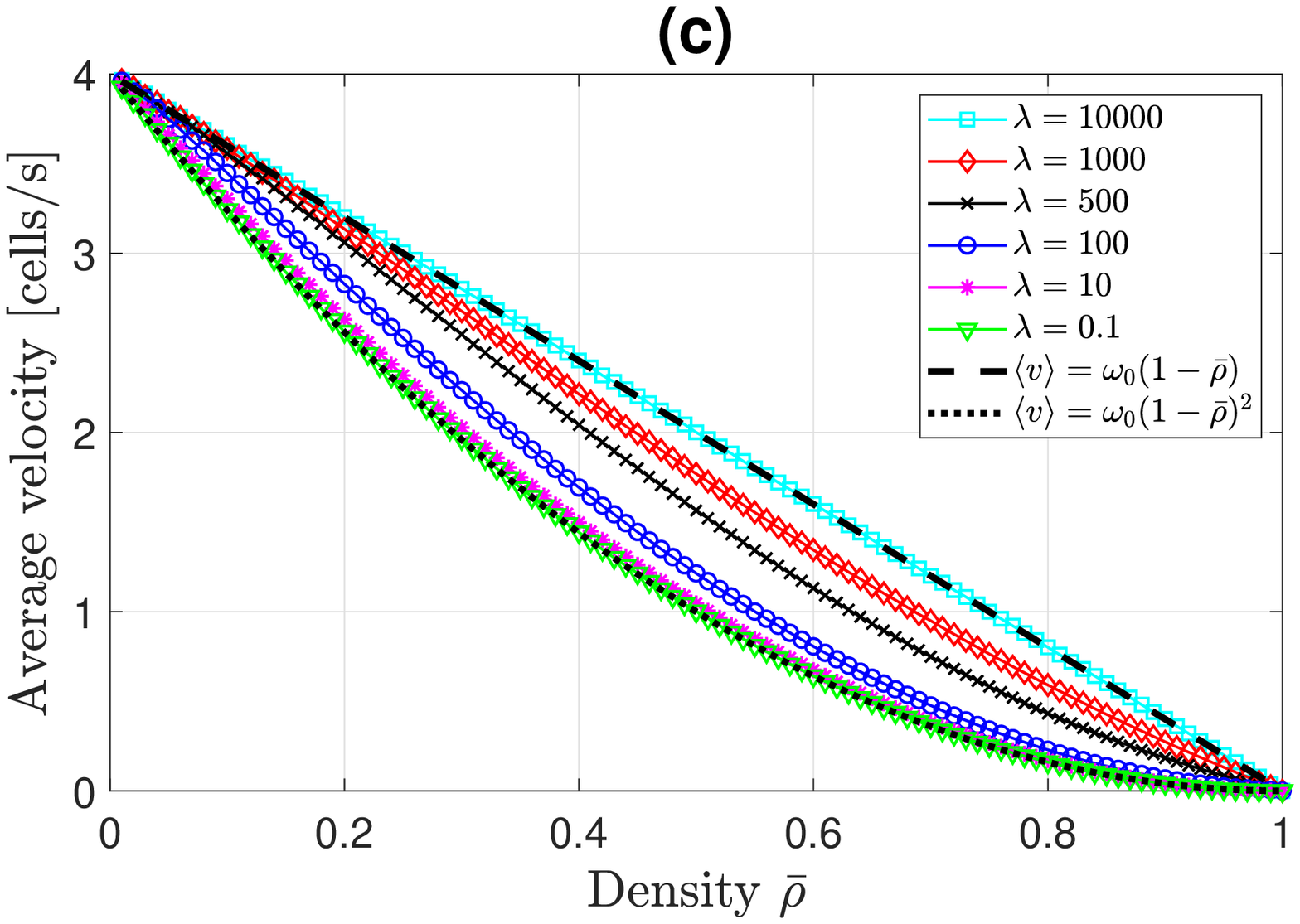} \hfill
\includegraphics[width=.46\textwidth]{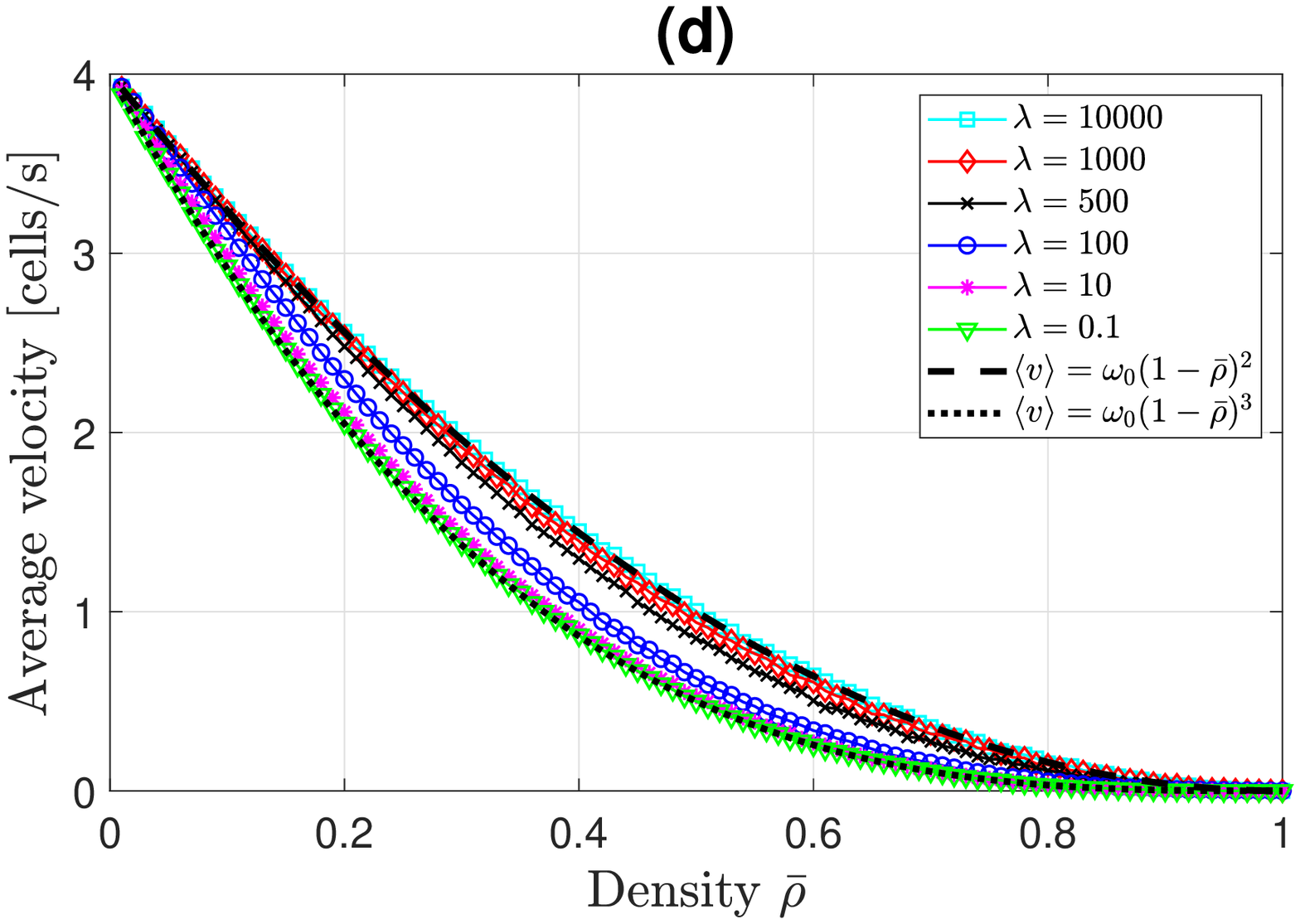}

\vspace{0.2cm}

\hspace{0.2cm} \includegraphics[width=.47\textwidth]{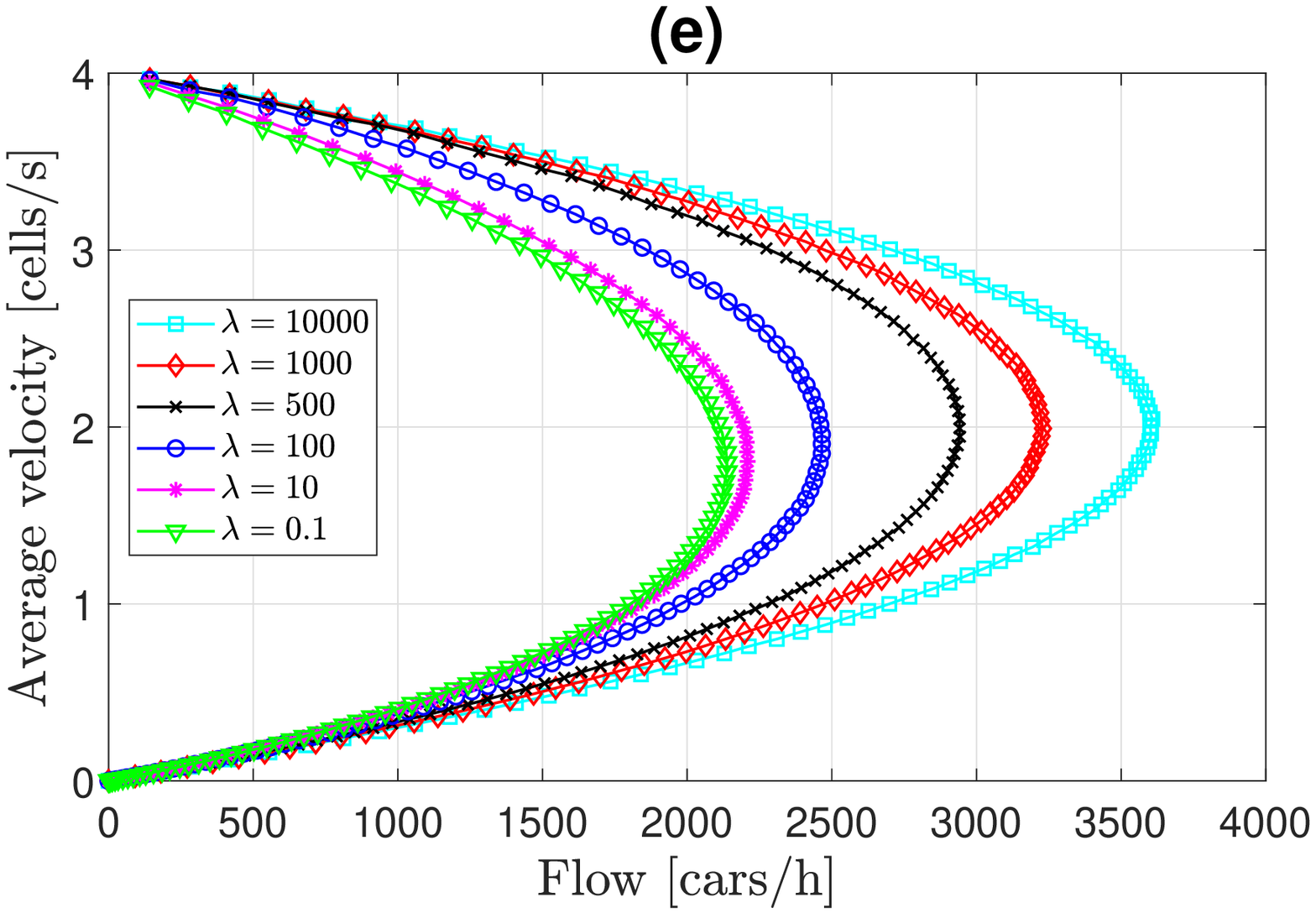} \hfill
\includegraphics[width=.47\textwidth]{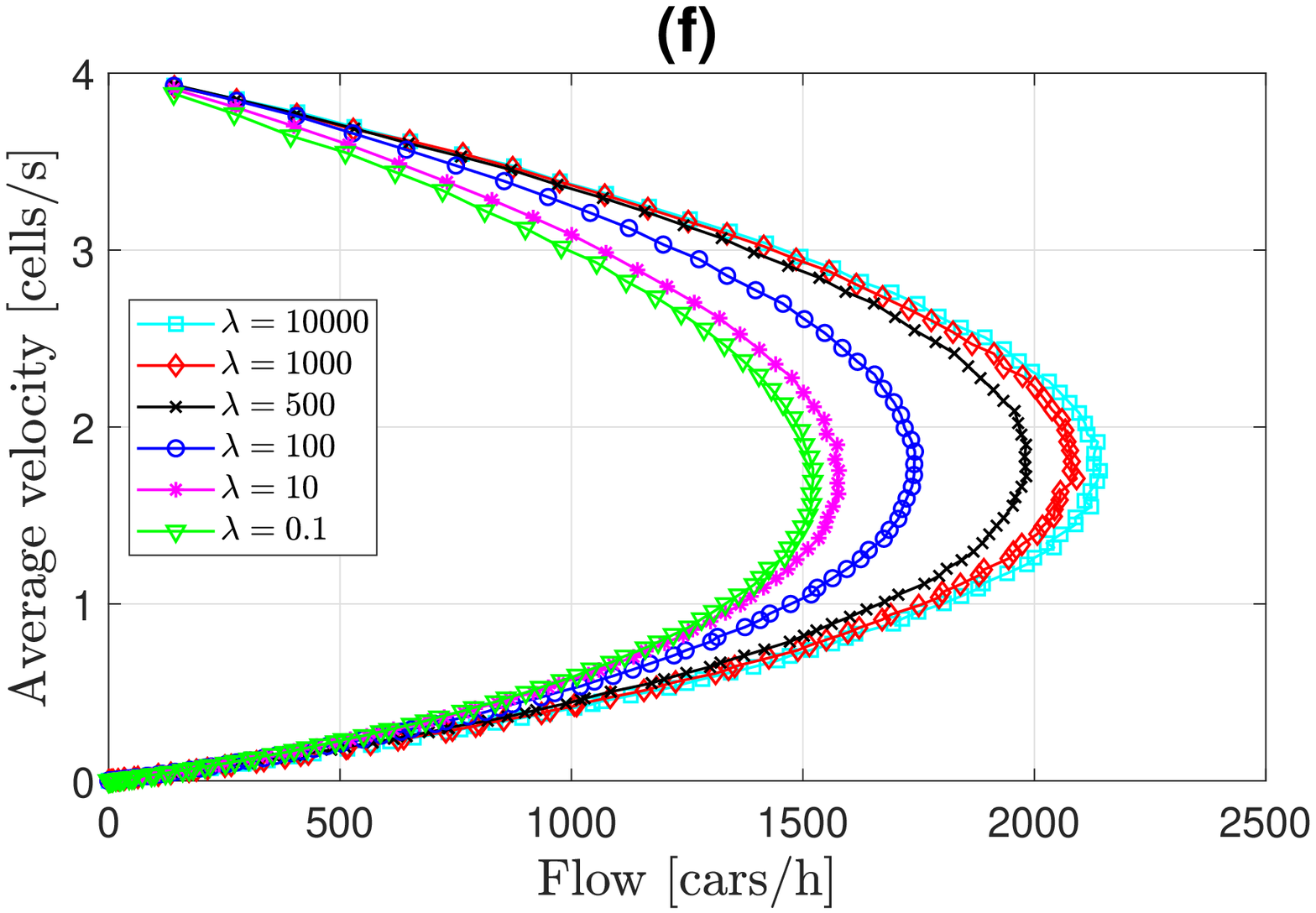}

\caption{Comparison results of the traffic flow on the one-lane highway with six different values of the interaction strength $\Eo$. In all KMC simulations, we take the highway distance of $\approx4.17$ miles ($\approx 6704$m, $\NL=1000$ cells), and the final time of $1$ h. (a)(b): Long-time averages of the density-flow relationship; (c)(d): Ensemble-averaged velocity of cars versus the density $\bar{\rho}$; (e)(f): Long-time averages of the flow-velocity relationship. Note that for $\Eo=0.1$ and $10^4$, the fluxes of the KMC simulations in (a) and (b) agree with the macroscopic averaged fluxes in \eqref{eq:fluxavg0} and \eqref{eq:fluxavginf} (shown as the dotted and dashed black curves, respectively). (left panels): Results of $\Jump=1$. (right panels): Results of $\Jump=2$. The slowdown  interaction function: $g(x) =1-x$ in \eqref{eq:g-1}.}
                                                           \label{fig:flowE}
\end{center}
\end{figure}

We use the accelerated KMC method to simulate the dynamics with different choices of $\lambda$.
In Fig.~\ref{fig:flowE}, the fundamental diagrams are plotted for different choices of $\Eo$:
\begin{equation}\label{eq:lambda}
  \Eo=0.1,\,\,10,\,\,100,\,\,500,\,\,1000,\,\,10000,
\end{equation}
and for different multiple move parameters $\Jump=1$ and $\Jump=2$.
All curves exhibit phase transitions between the free-flow phase and the jammed phase.

\begin{figure}[p]
\begin{center}
\includegraphics[width=.48\textwidth]{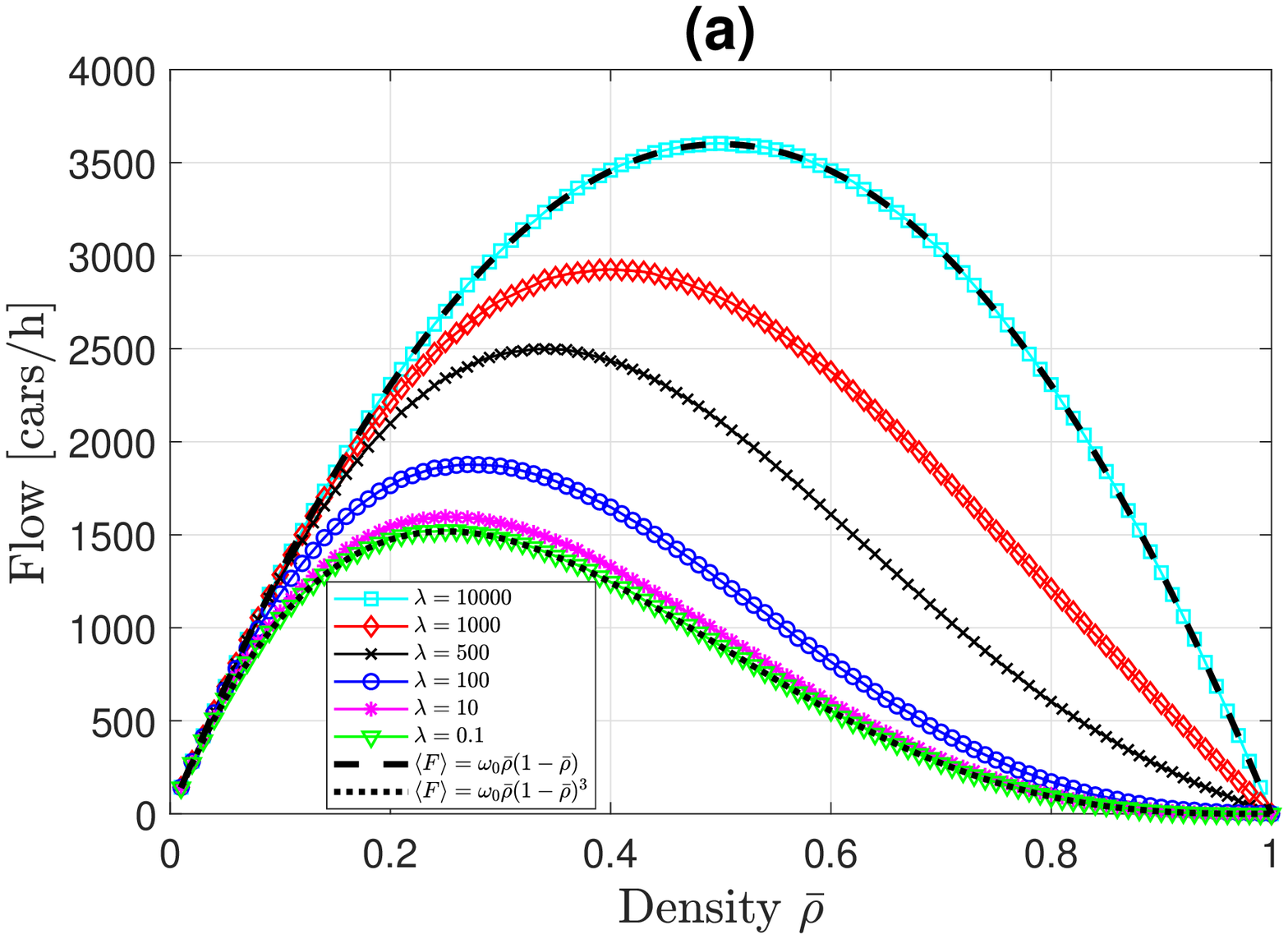} \hfill
\includegraphics[width=.48\textwidth]{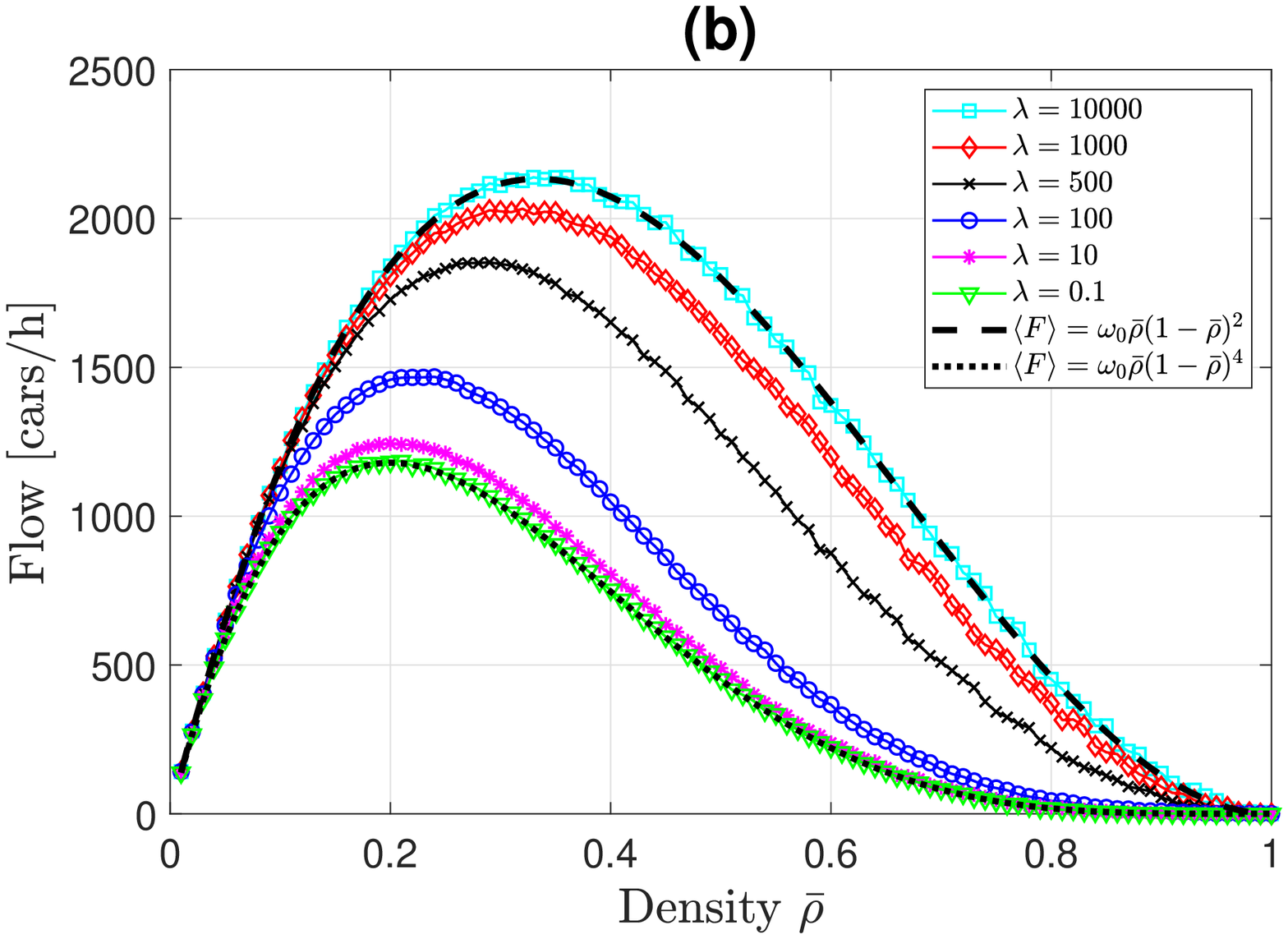}

\vspace{0.2cm}

\hspace{0.2cm} \includegraphics[width=.46\textwidth]{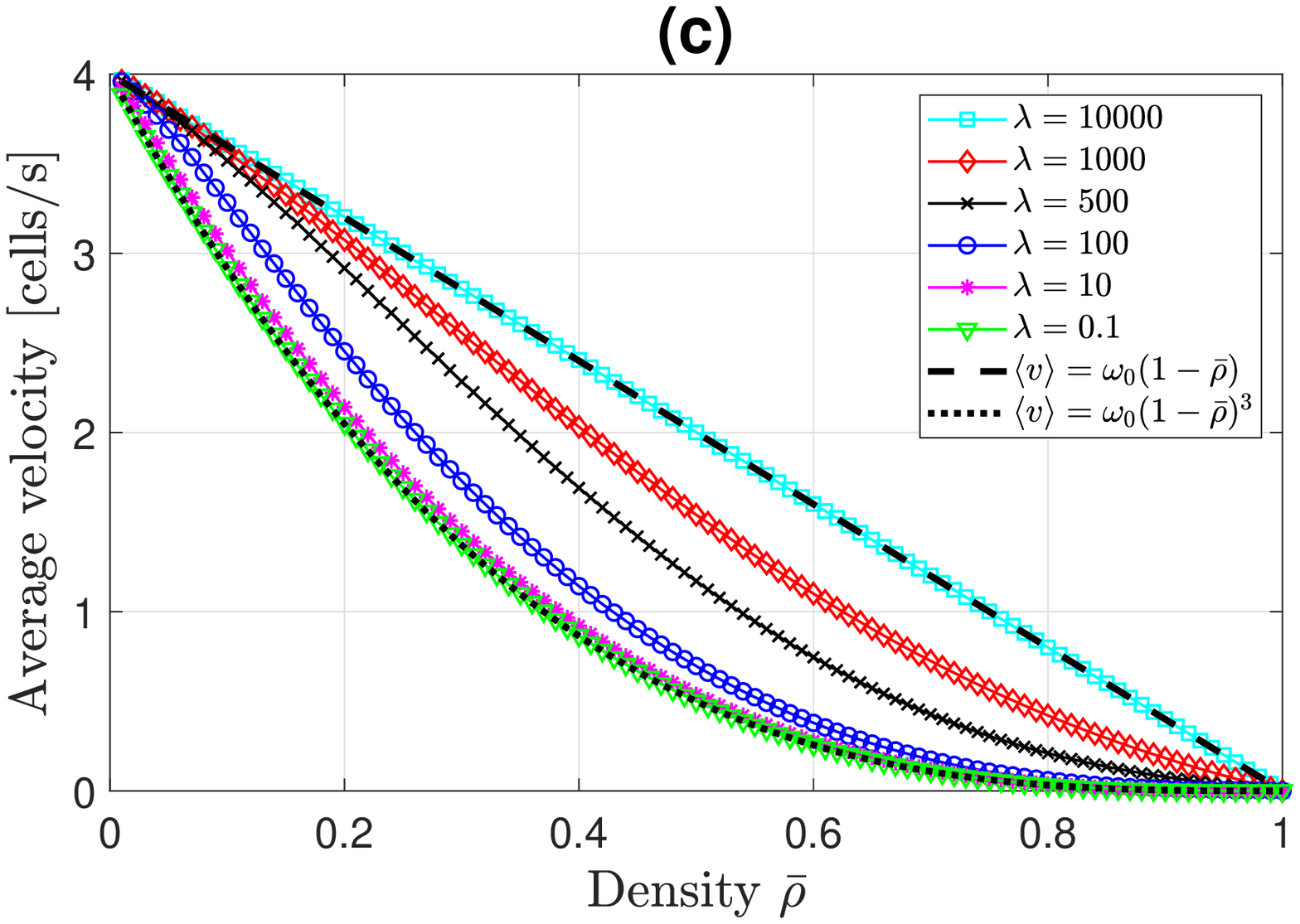} \hfill
\includegraphics[width=.46\textwidth]{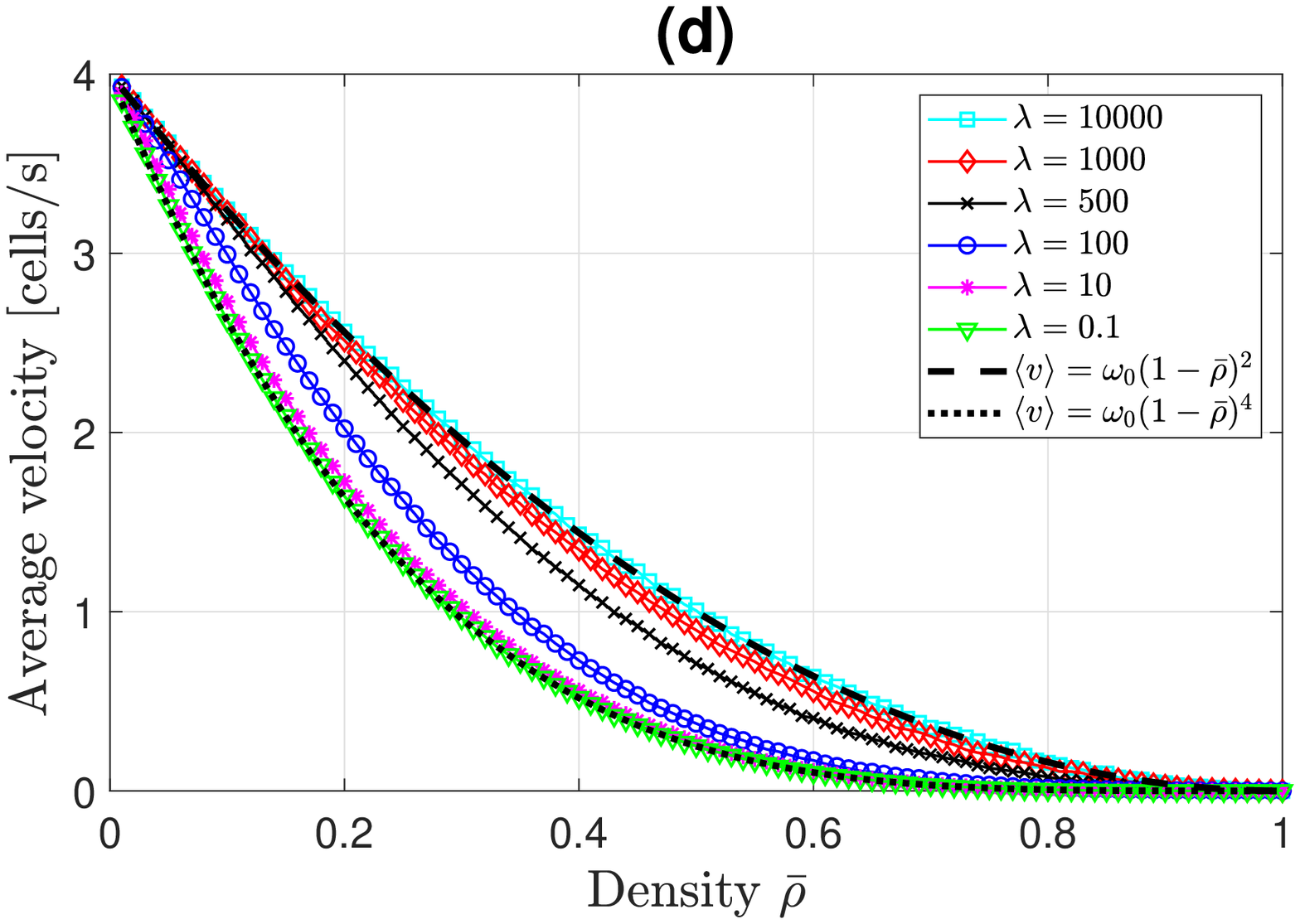}

\vspace{0.2cm}

\hspace{0.2cm} \includegraphics[width=.47\textwidth]{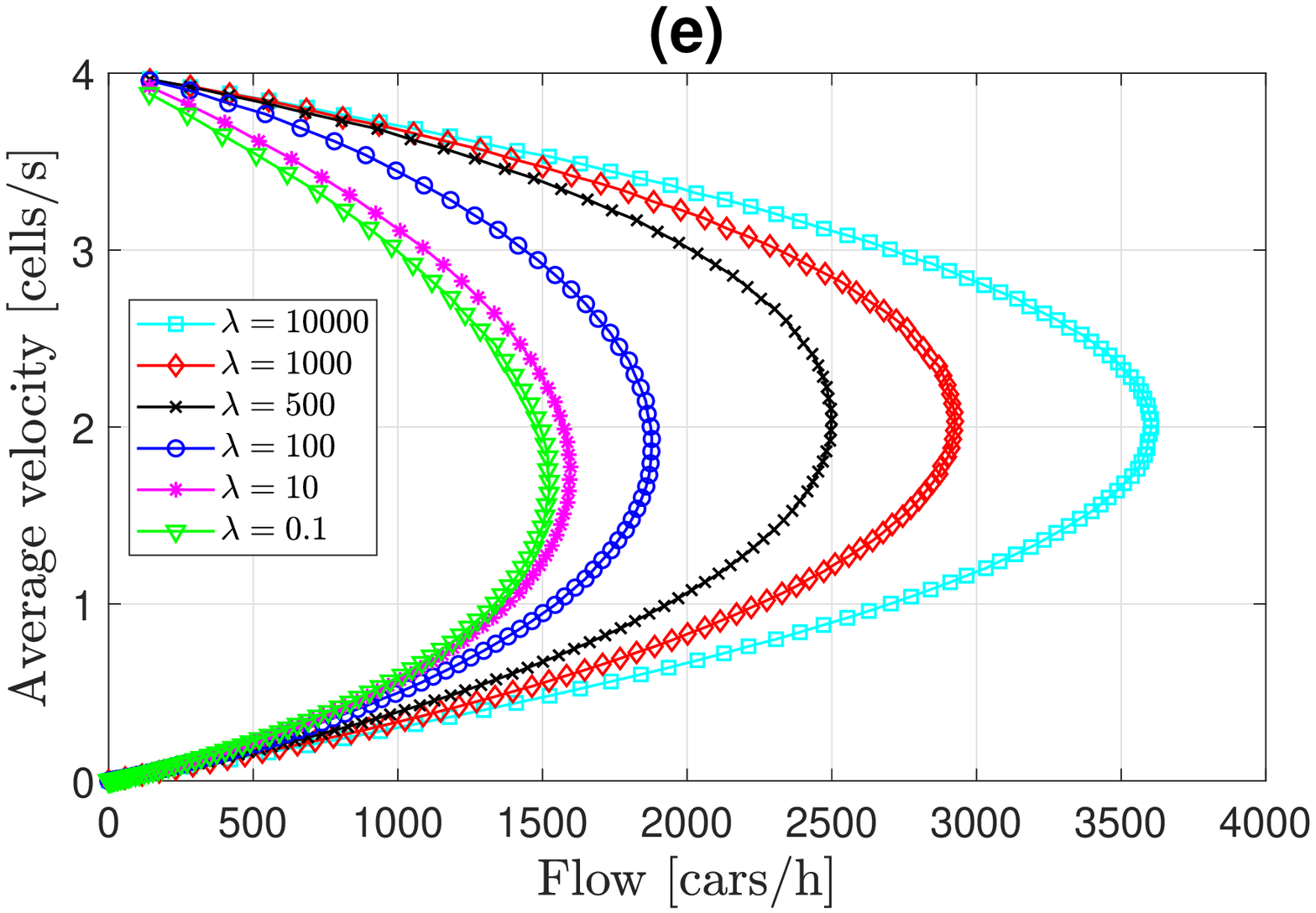} \hfill
\includegraphics[width=.47\textwidth]{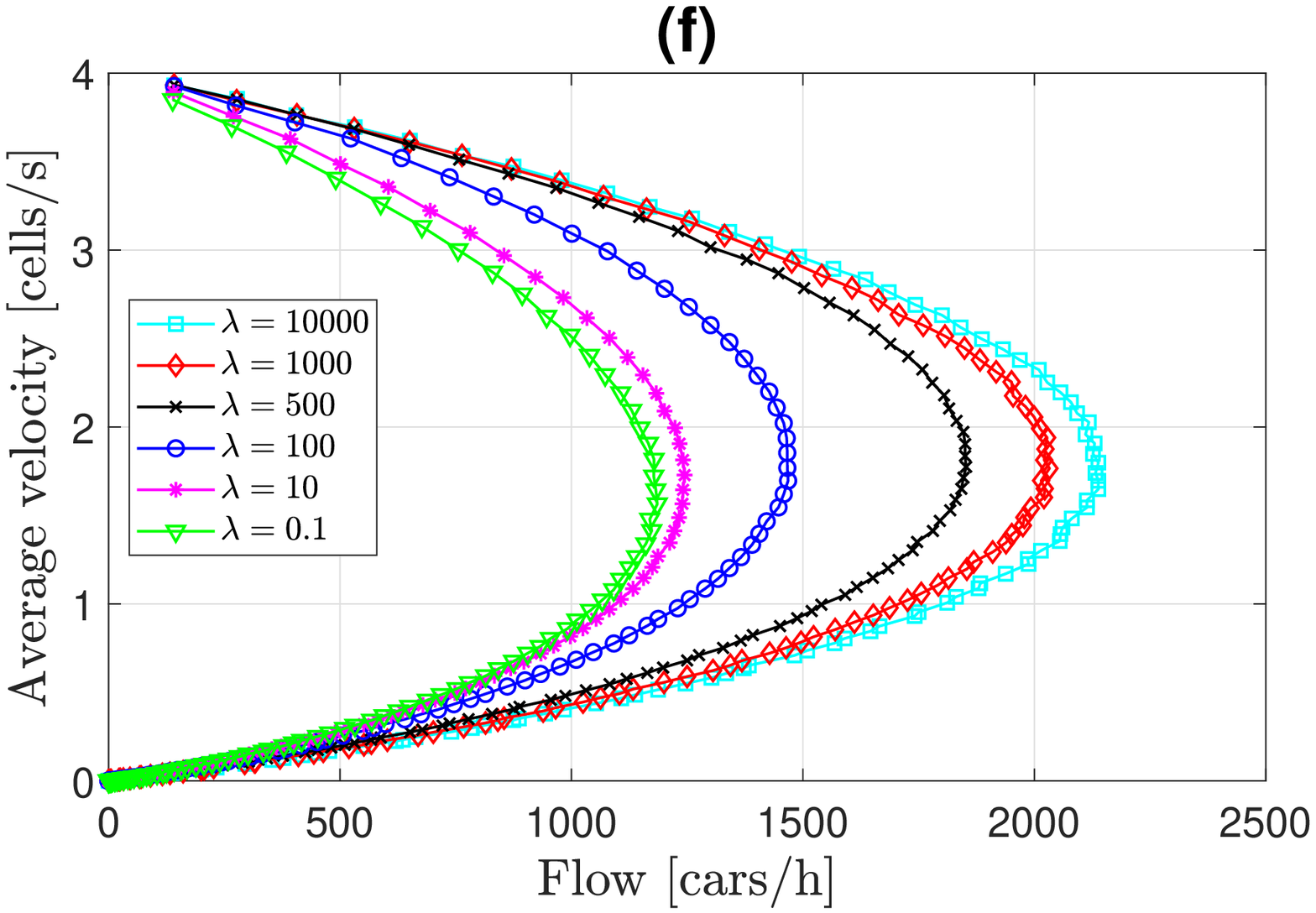}

\caption{Comparison results of the traffic flow on the one-lane highway with six different values of the interaction strength $\Eo$. (a)(b): Long-time averages of the density-flow relationship; (c)(d): Ensemble-averaged velocity of cars versus the density $\bar{\rho}$; (e)(f): Long-time averages of the flow-velocity relationship. Note that for $\Eo=0.1$ and $10^4$, the fluxes of the KMC simulations in (a) and (b) agree with the macroscopic averaged fluxes in \eqref{eq:genavg0} and \eqref{eq:fluxavginf} (shown as the dotted and dashed black curves, respectively). (left panels): Results of $\Jump=1$. (right panels): Results of $\Jump=2$. All parameters are the same as in Fig.~\ref{fig:flowE} except that the slowdown interaction function: $g(x) =(1-x)^2$ in \eqref{eq:g-2}.}
                                                           \label{fig:flowE2}
\end{center}
\end{figure}

In Figs.~\ref{fig:flowE}(a) and (b), we plot the fundamental diagrams on the averaged fluxes $\Fbar$ against the averaged density $\bar\rho$ with different $\Eo$ in \eqref{eq:lambda} and $\Jump=1, 2$, respectively. They all share certain characteristics: a nearly linear increase of the flow at low averaged densities (which corresponds to the free-flow regime), a single maximum of the flow reached at a critical density $\rhoc^\Eo$, and a right-skewed asymmetry (namely $\rhoc^\Eo<1/2$).
We also observe that both the value of the critical density $\rhoc^\Eo$ and the maximum value of the flow $\Fbar$ tend to decrease with decreasing $\Eo$.
In particular, for the cases of $\Eo=0.1$ and $10^4$, the fluxes of the KMC simulations agree with the limiting fluxes in \eqref{eq:fluxavg0} and \eqref{eq:fluxavginf}, (shown as the dotted and dashed black curves, respectively).

In Fig.~\ref{fig:flowE}(a) for the case of $\Jump=1$, the density-flow curves for the cases of $\Eo=0.1$ and $10^4$ take their maxima at the critical density around $\rhoc=\frac13$ and $\frac12$, respectively, which is consistent with the limiting cases \eqref{eq:fluxavg0} and \eqref{eq:fluxavginf}
\begin{equation}\label{eq:rhocrit}
  \lim_{\Eo\to0}\rho_c^\Eo=\frac{1}{\Jump+2},\qquad
  \lim_{\Eo\to\infty}\rho_c^\Eo=\frac{1}{\Jump+1}.
\end{equation}
The maximum flux for $\Eo=0.1$ at $\rhoc=\frac13$ is about $3600 \omega_0 \cdot\frac{4}{27} \approx 2133$ cars per hour (recall that $\omega_0=4$s$^{-1}$),
and the maximum flux for $\Eo=10^4$ at $\rhoc=\frac12$ is about $3600 \omega_0 \cdot\frac{1}{4} = 3600$ cars per hour.
Similarly, in Fig.~\ref{fig:flowE}(b) for the case of $\Jump=2$, the density-flow curves for the cases of $\Eo=0.1$ and $10^4$ take their maxima $3600 \omega_0 \cdot\frac{27}{256} \approx 1519$ and $3600 \omega_0 \cdot\frac{4}{27} \approx 2133$ at the critical density $\rhoc=\frac14$ and $\frac13$, respectively.

Figs.~\ref{fig:flowE}(c) and (d) show the fundamental diagrams of the density-velocity relationship for $\Jump=1$ and $\Jump=2$, respectively. For the special case of $\Jump=1$ and $\Eo=10^4$ in Fig.~\ref{fig:flowE}(c), the ensemble-averaged velocity $\vbar$ decreases linearly as the averaged density $\bar\rho$ increases, which is consistent with the classical  Lighthill-Whitham-Richards (LWR) model  \cite{LiW55, Whi74}. All other cases in both Figs.~\ref{fig:flowE}(c) and (d) show that in the free-flow regime the ensemble-averaged velocity $\vbar$ decreases approximately linearly from the full speed of $4$ cells per second ($\approx 26.8$ m/s or $60$ miles/h) as $\bar{\rho}$ increases and the chance of interaction between cars gets higher. As $\Eo$ decreases, when $\bar{\rho}$ is larger than the critical point $\rhoc$, the average velocity drops down to zero and the density-velocity curve is convex. This linear relationship follows the Greenshields model \cite{Gre35} and the convex relationship belongs to the Underwood model \cite{Und61}.  For the cases of $\Eo=0.1$ and $10^4$  in both Figs.~\ref{fig:flowE}(c) and (d), the ensemble-averaged velocity of the KMC simulations agrees with limits \eqref{eq:fluxavg0} and \eqref{eq:fluxavginf}, (shown as the dotted and dashed black curves, respectively).

Figs.~\ref{fig:flowE}(e) and (f) show the fundamental diagrams of the flow-velocity relationship for $\Jump=1$ and $\Jump=2$, respectively, which plot the ensemble-averaged velocity $\vbar$ versus the averaged flow $\Fbar$.
In Fig.~\ref{fig:flowE}(e) of $\Jump=1$, for the case of $\Eo=10^4$ (shown as cyan squares), the flow $\Fbar$ reaches its maximum $\approx 3600$ cars per hour when the ensemble-averaged velocity $\vbar$ is at a critical value $\vbar_{\rm crit} \approx 2.0$ cells per second ($\approx 13.4$ m/s or $30$ miles/h). As $\Eo$ decreases to $0.1$, the maximum value of the flow decreases a lot to $\approx 2133$ cars per hour and the critical value $\vbar_{\rm crit}$ decreases a bit and becomes slightly lower than $2$ cells per second. For the case of $\Jump=2$ in Fig.~\ref{fig:flowE}(f), when $\Eo$ decreases, the maximum flow $\Fbar$ decreases from $\approx 2133$ cars per hour for the case of $\Eo=10^4$ down to $\approx 1519$ cars per hour for the case of $\Eo=0.1$. The critical value $\vbar_{\rm crit}$ where the flow $\Fbar$ reaches its maximum do not change too much. The results compare favorably with observed data in \cite{Wie95}.

\subsection{General slowdown relations}
In the following, we perform numerical experiments on a family of CA models with the same kernel \eqref{eq:kiBS} but a different slowdown relation
\begin{equation}\label{eq:g-2}
g(x) = \begin{cases}(1-x)^2& x\in[0,1],\\0&\text{otherwise}.\end{cases}
\end{equation}
We apply the accelerated KMC method and generate fundamental diagrams of the dynamics with $\Eo$ in \eqref{eq:lambda} and $\Jump=1, 2$.

Let us summarize the expected behaviors of the two extreme cases.
When $\Eo\to0$, we have
\begin{equation}\label{eq:genavg0}
  \Fbar(\rhob)=\omega_0\rhob(1-\rhob)^{\Jump+2},
  \qquad \vbar(\rhob)=\omega_0(1-\rhob)^{\Jump+2},
  \qquad \lim_{\Eo\to0}\rhoc^\Eo=\frac{1}{\Jump+3}.
\end{equation}
When $\Eo\to\infty$, $\Fbar$ and $\vbar$ should behave the same as in \eqref{eq:fluxavginf}.

In Fig.~\ref{fig:flowE2}, we observe the essentially same phenomena about the fundamental diagrams as shown in Fig.~\ref{fig:flowE}. In particular, the averaged results for the case of $\Eo=10^4$ (shown as cyan squares) in all panels of Fig.~\ref{fig:flowE2} are statistically equal to the ones shown in Fig.~\ref{fig:flowE}. But for $\Eo=0.1$ (shown as green ``$\triangledown$" signs), the density-flow curves in Fig.~\ref{fig:flowE2}(a)(b) for the cases of $\Jump=1$ and $\Jump=2$ take their maxima at the critical density $\rhoc=\frac14$ and $\frac15$, respectively, which is consistent with \eqref{eq:genavg0}.
Moreover, for all cases of $\Eo < 10^4$, the density-flow curves shown in Fig.~\ref{fig:flowE2}(a)(b) and the density-velocity curves in Fig.~\ref{fig:flowE2}(c)(d) are below the corresponding curves shown in Fig.~\ref{fig:flowE}(a)-(d), respectively. Meanwhile, the flow-velocity curves in Fig.~\ref{fig:flowE2}(e)(f) are left of the corresponding ones shown in Fig.~\ref{fig:flowE}(e)(f). This is because the function $g(x)=(1-x)^2$ introduces a stronger slow-down effect than $g(x)=1-x$ does such that for the same value of $\rhob$, both the averaged flux and the ensemble-averaged velocity are reduced.

\subsection{Different multiple move parameters}
\label{subsect.jump}

\begin{figure}[p]
\begin{center}
\includegraphics[width=.48\textwidth]{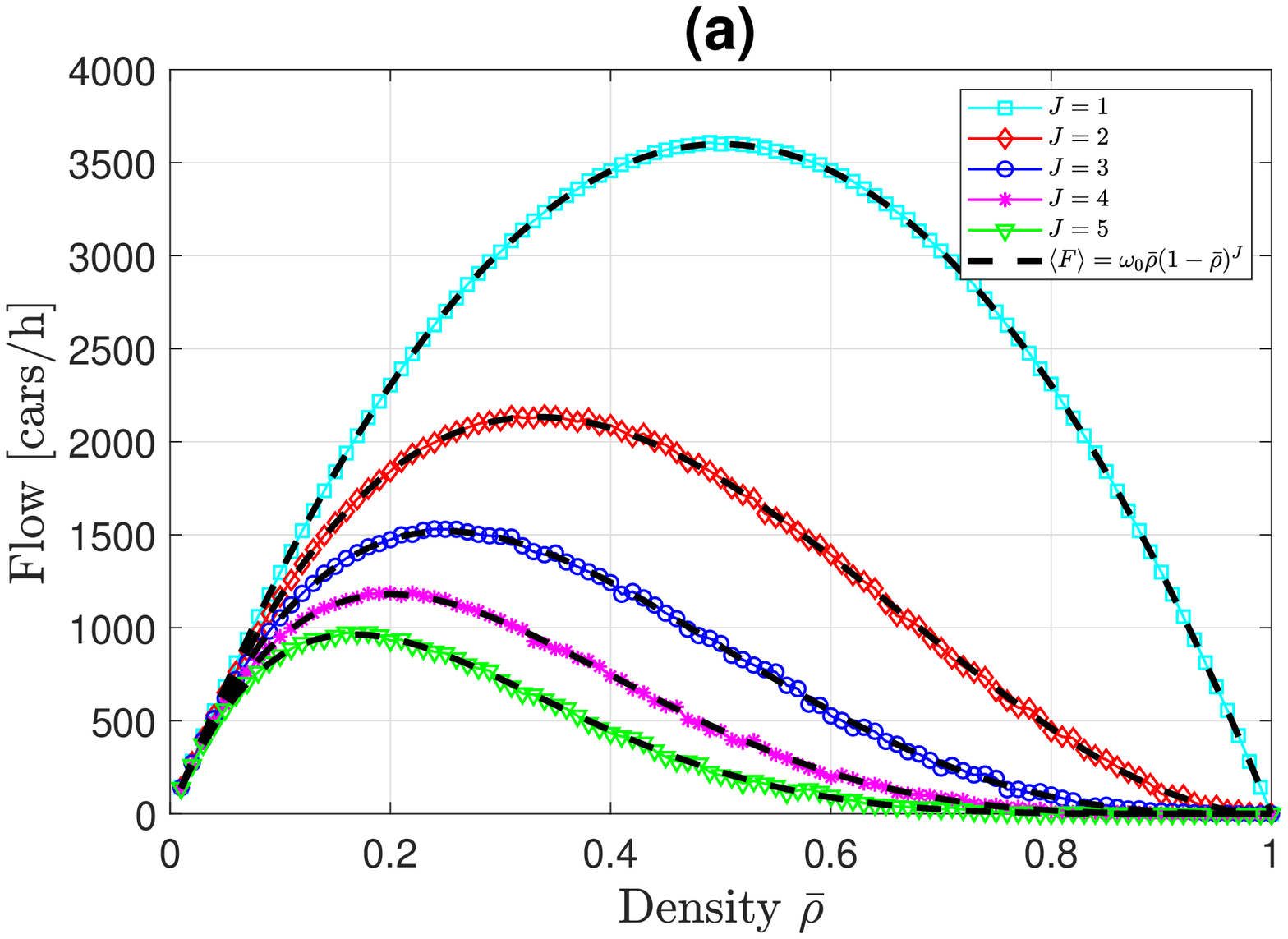} \hfill
\includegraphics[width=.48\textwidth]{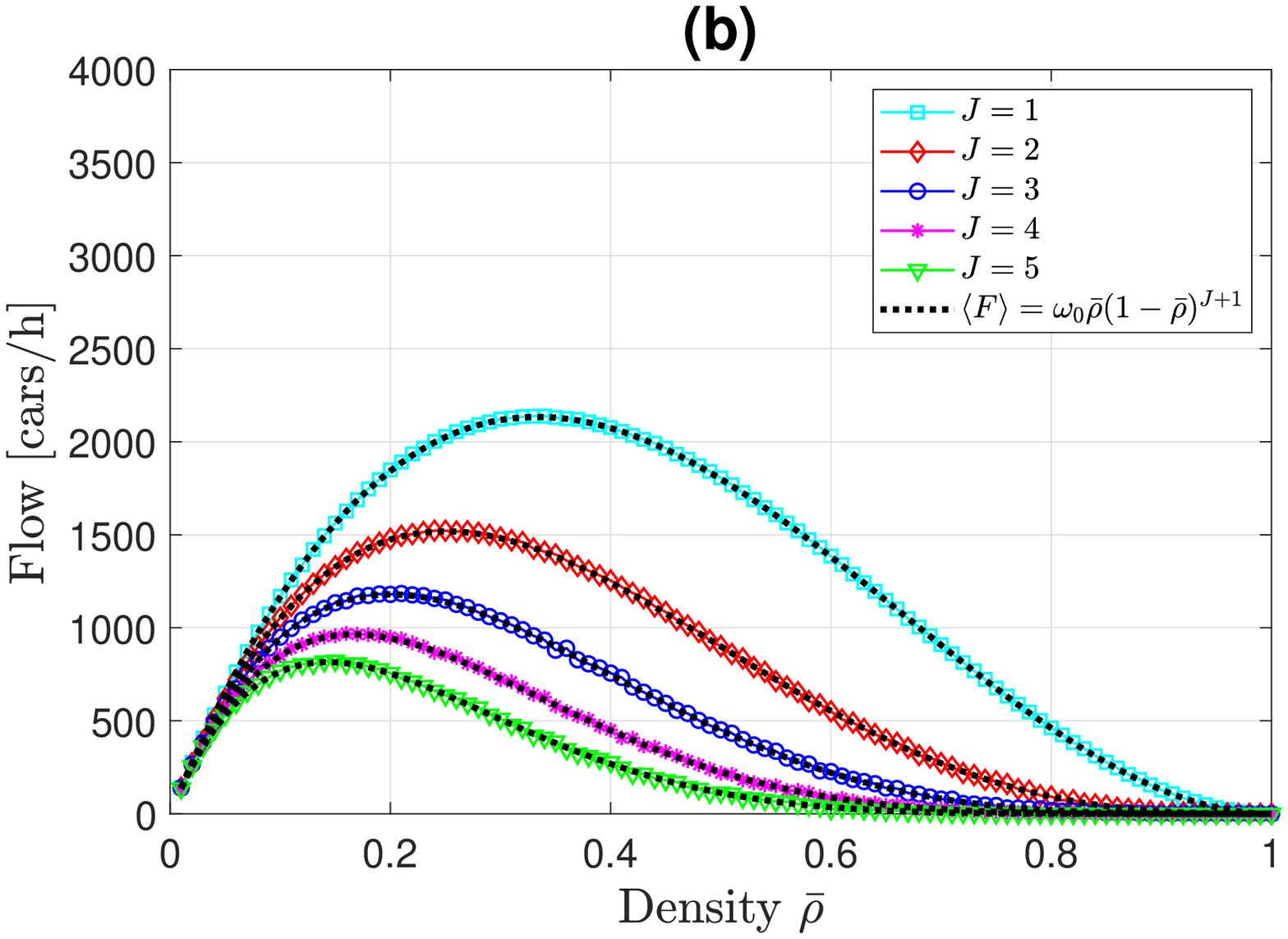}

\vspace{0.2cm}

\hspace{0.2cm} \includegraphics[width=.46\textwidth]{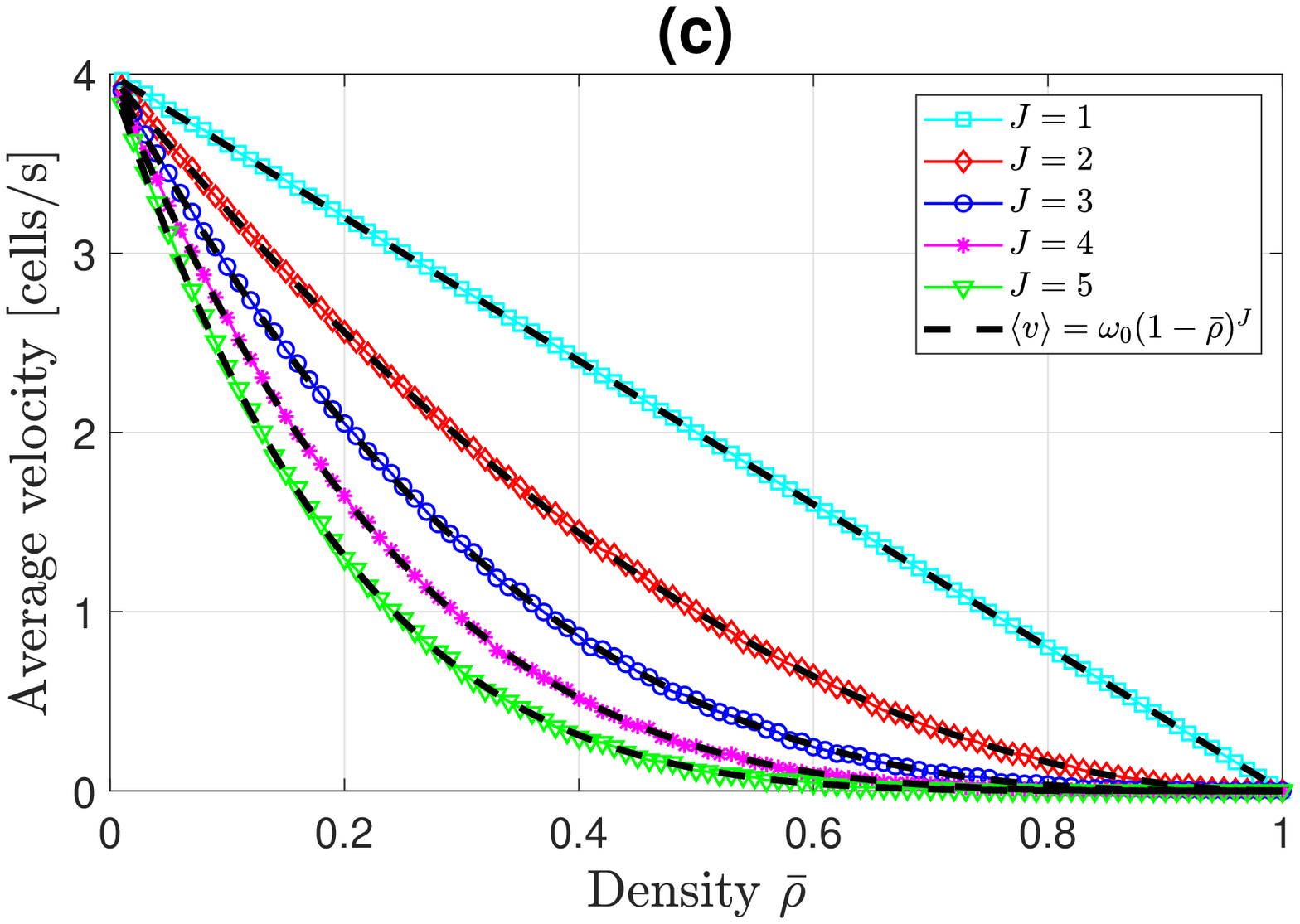} \hfill
\includegraphics[width=.46\textwidth]{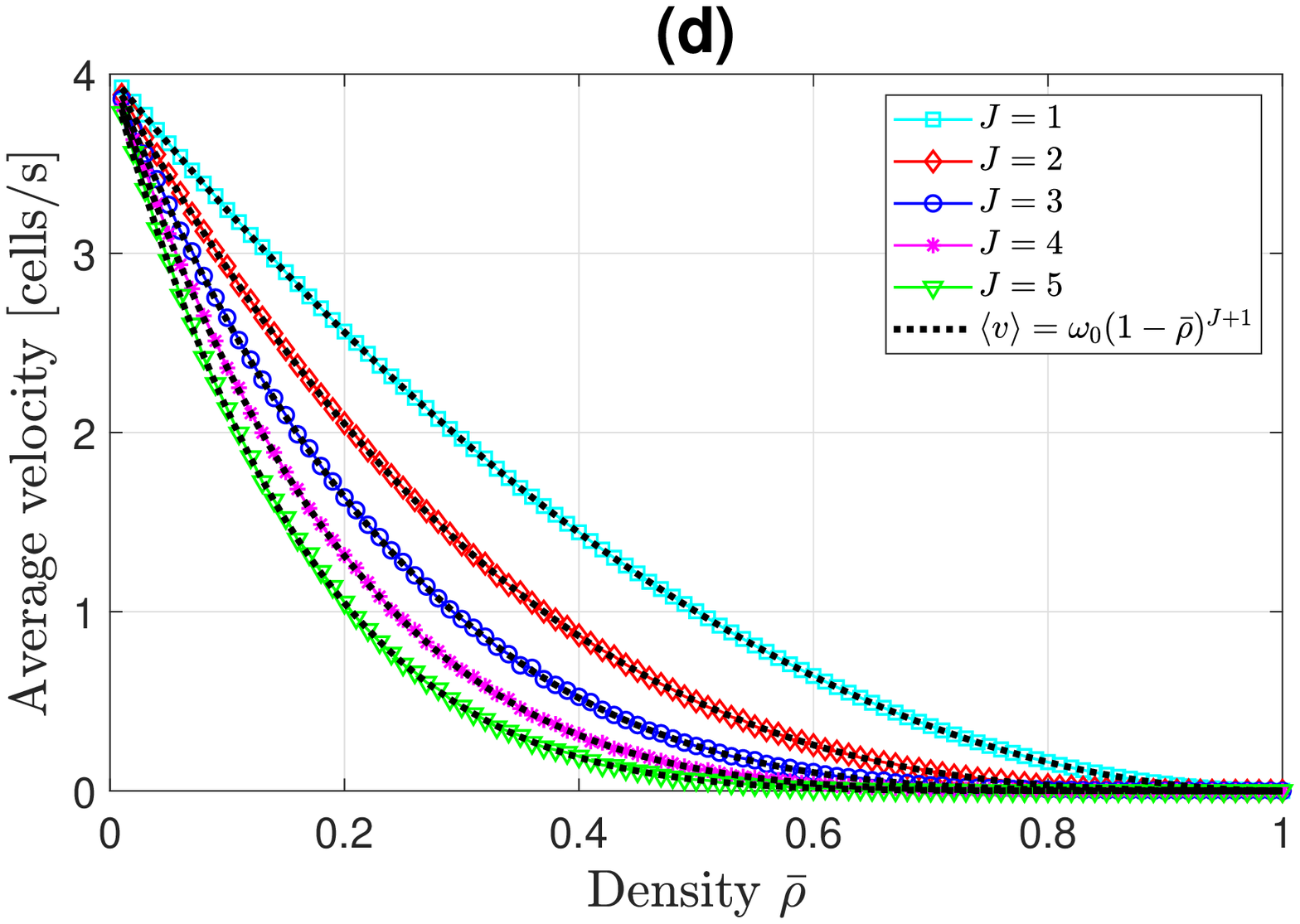}

\vspace{0.2cm}

\hspace{0.2cm} \includegraphics[width=.47\textwidth]{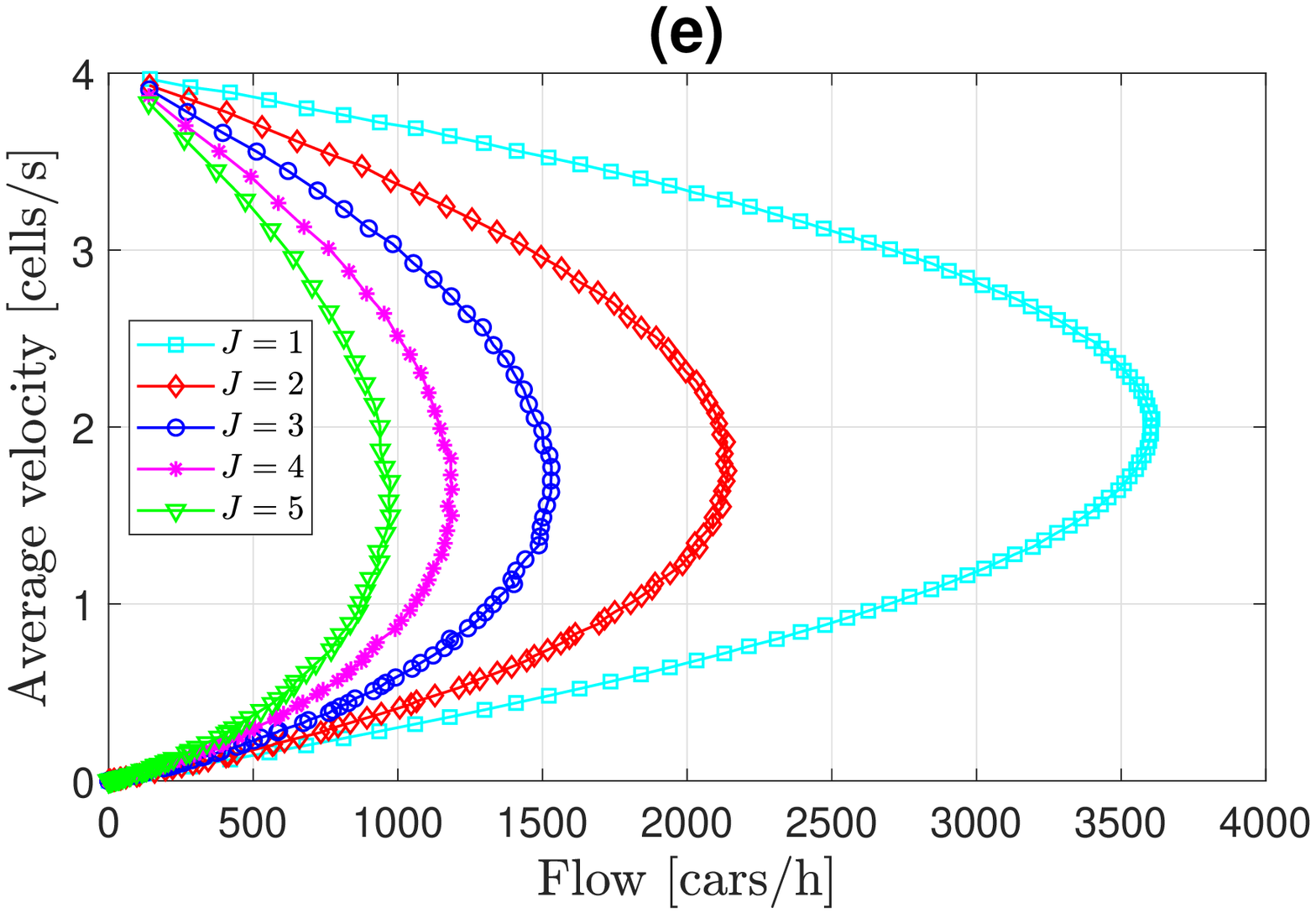} \hfill
\includegraphics[width=.47\textwidth]{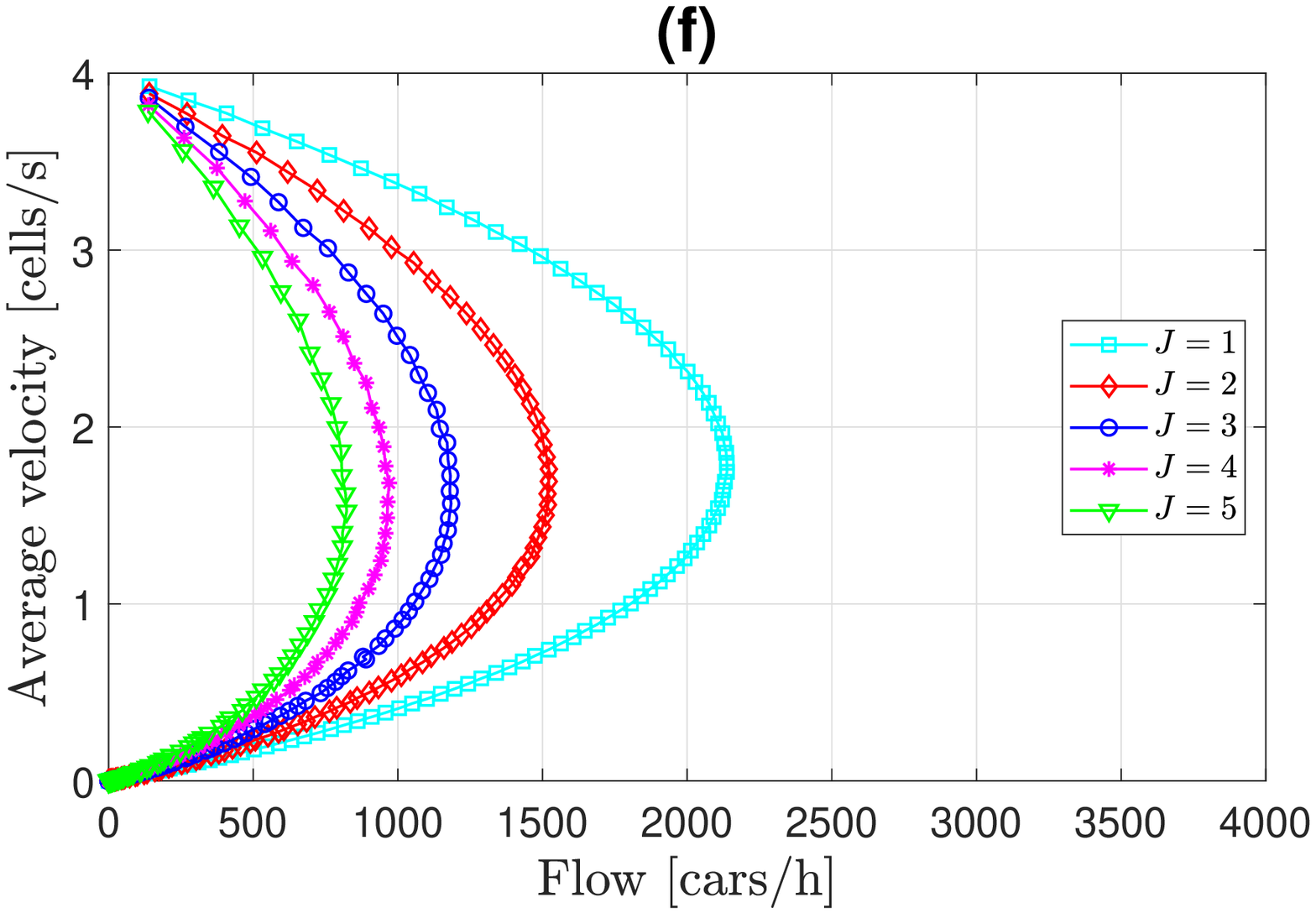}

\caption{Comparison results of the traffic flow on the one-lane highway with five different values of the multiple move parameter $\Jump$. (a)(b): Long-time averages of the density-flow relationship; (c)(d): Ensemble-averaged velocity of cars versus the density $\bar{\rho}$; (e)(f): Long-time averages of the flow-velocity relationship. (left panels): Results of $\Eo=10^4$. (right panels): Results of $\Eo=0.1$. Note that for each value of $\Jump=1,2,\ldots, 5$, the fluxes of the KMC simulations in (a) and (b) agree with the macroscopic averaged fluxes in \eqref{eq:fluxavginf} and \eqref{eq:fluxavg0} (shown as the dashed and dotted black curves, respectively). The slowdown  interaction function: $g(x) =1-x$ in \eqref{eq:g-1}.}
                                                           \label{fig:flowJ}
\end{center}
\end{figure}

Finally, we show the effects of the multiple move parameter $\Jump$ on the flows in more detail in Fig.~\ref{fig:flowJ}. We first take the slowdown relation $g(x)=1-x$ in \eqref{eq:g-1} and compare the results of $\Eo=10^4$ and $\Eo=0.1$, respectively.

\begin{figure}[p]
\begin{center}
\includegraphics[width=.48\textwidth]{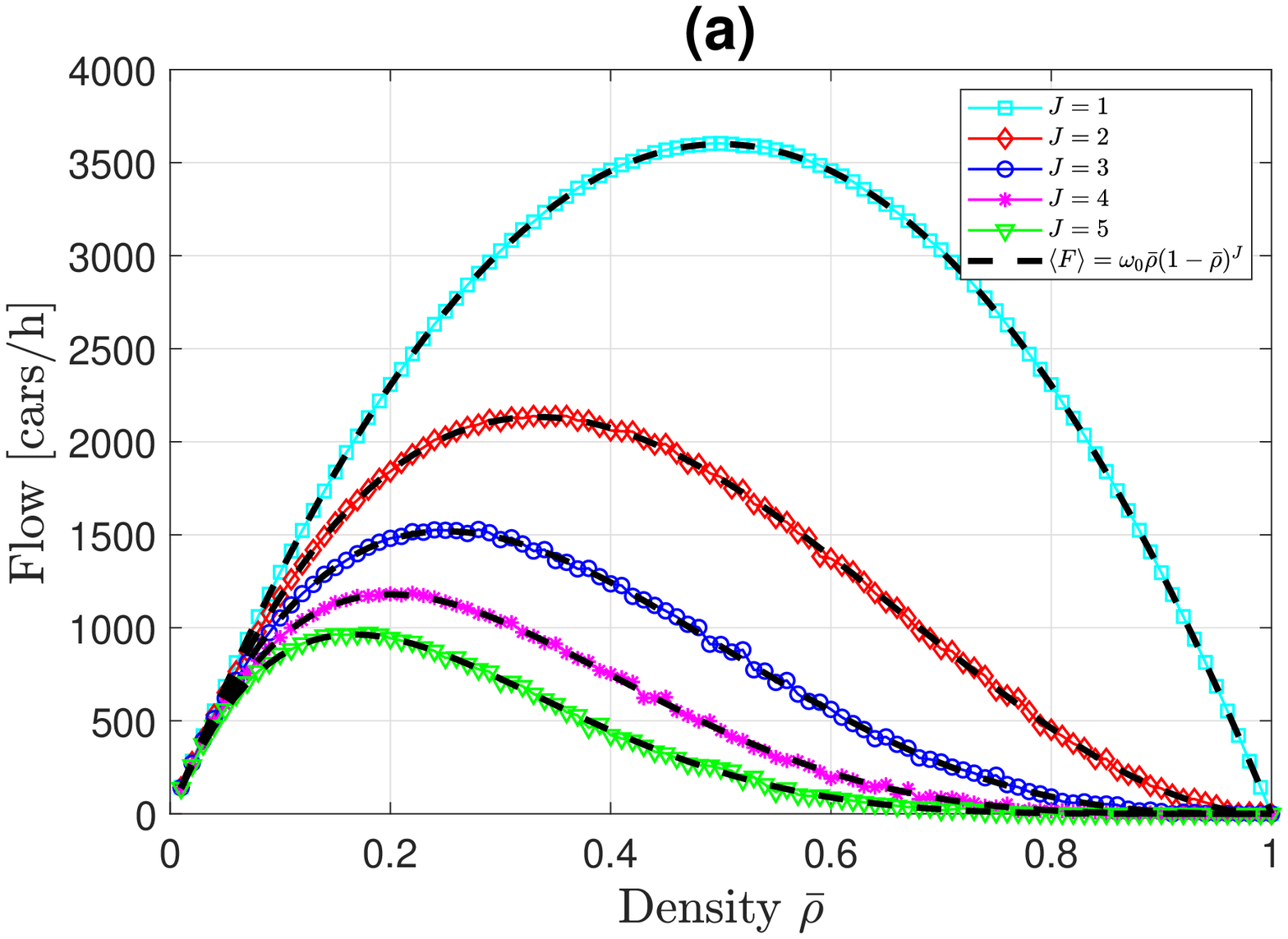} \hfill
\includegraphics[width=.48\textwidth]{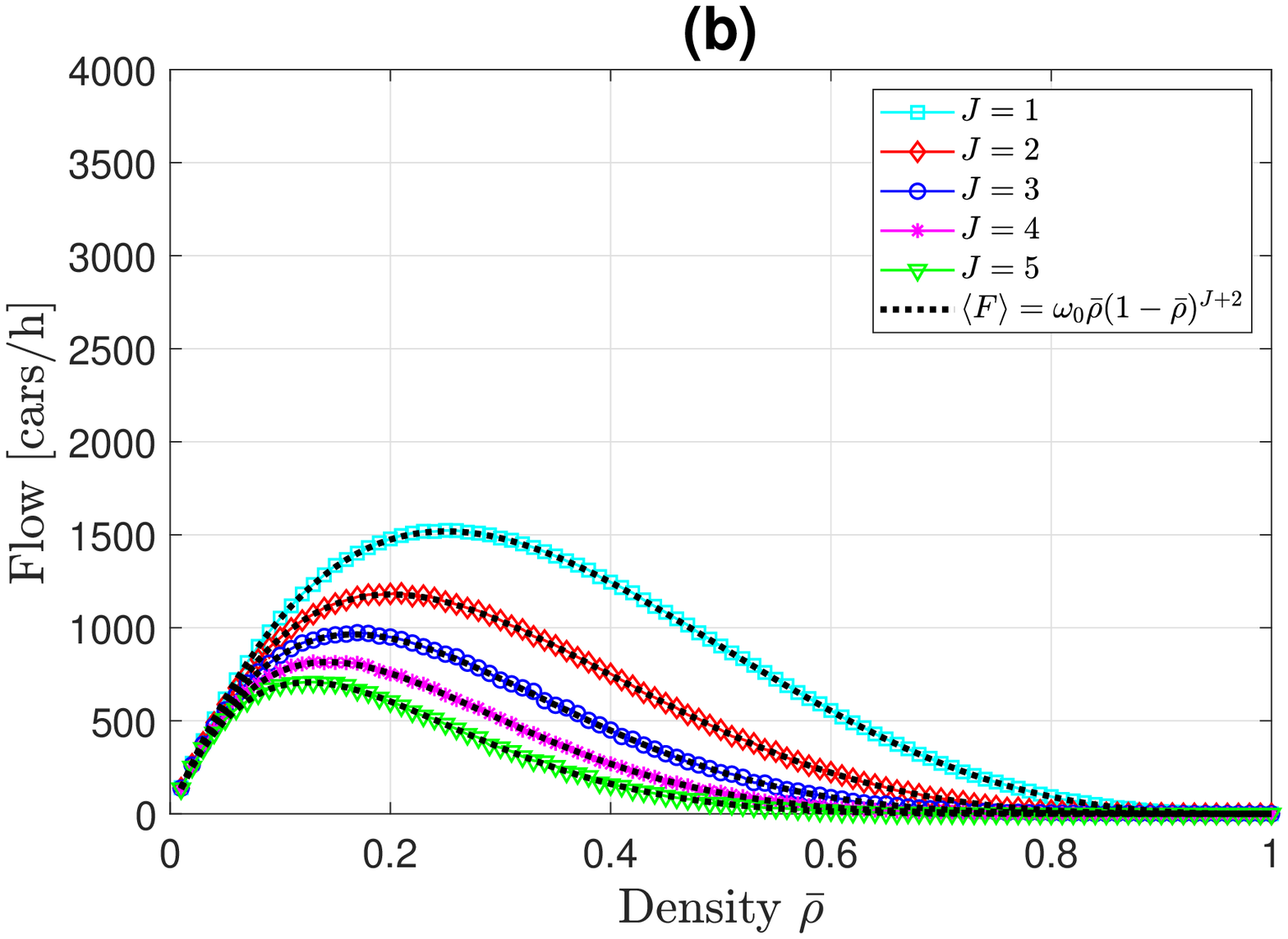}

\vspace{0.2cm}

\hspace{0.2cm} \includegraphics[width=.46\textwidth]{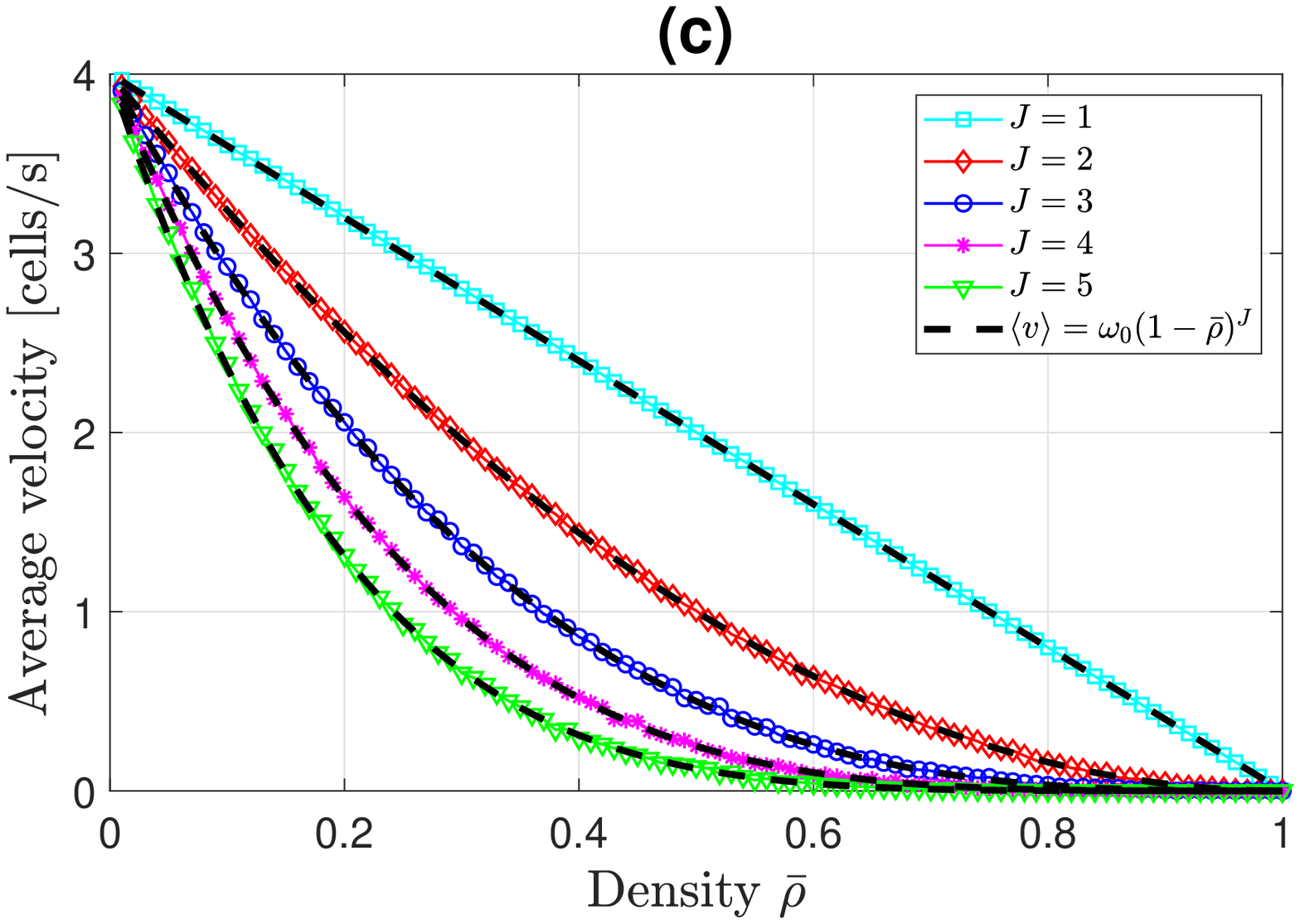} \hfill
\includegraphics[width=.46\textwidth]{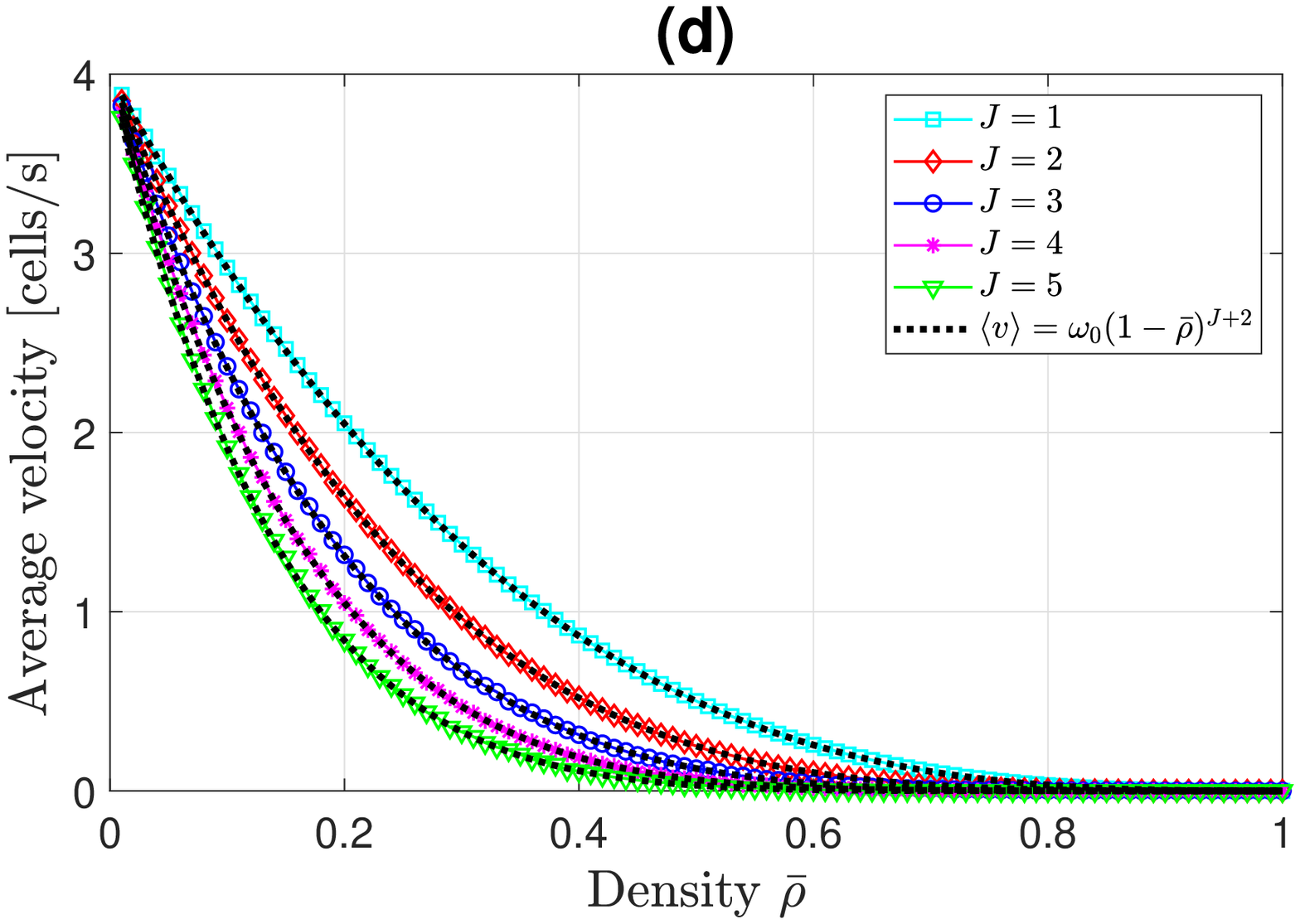}

\vspace{0.2cm}

\hspace{0.2cm} \includegraphics[width=.47\textwidth]{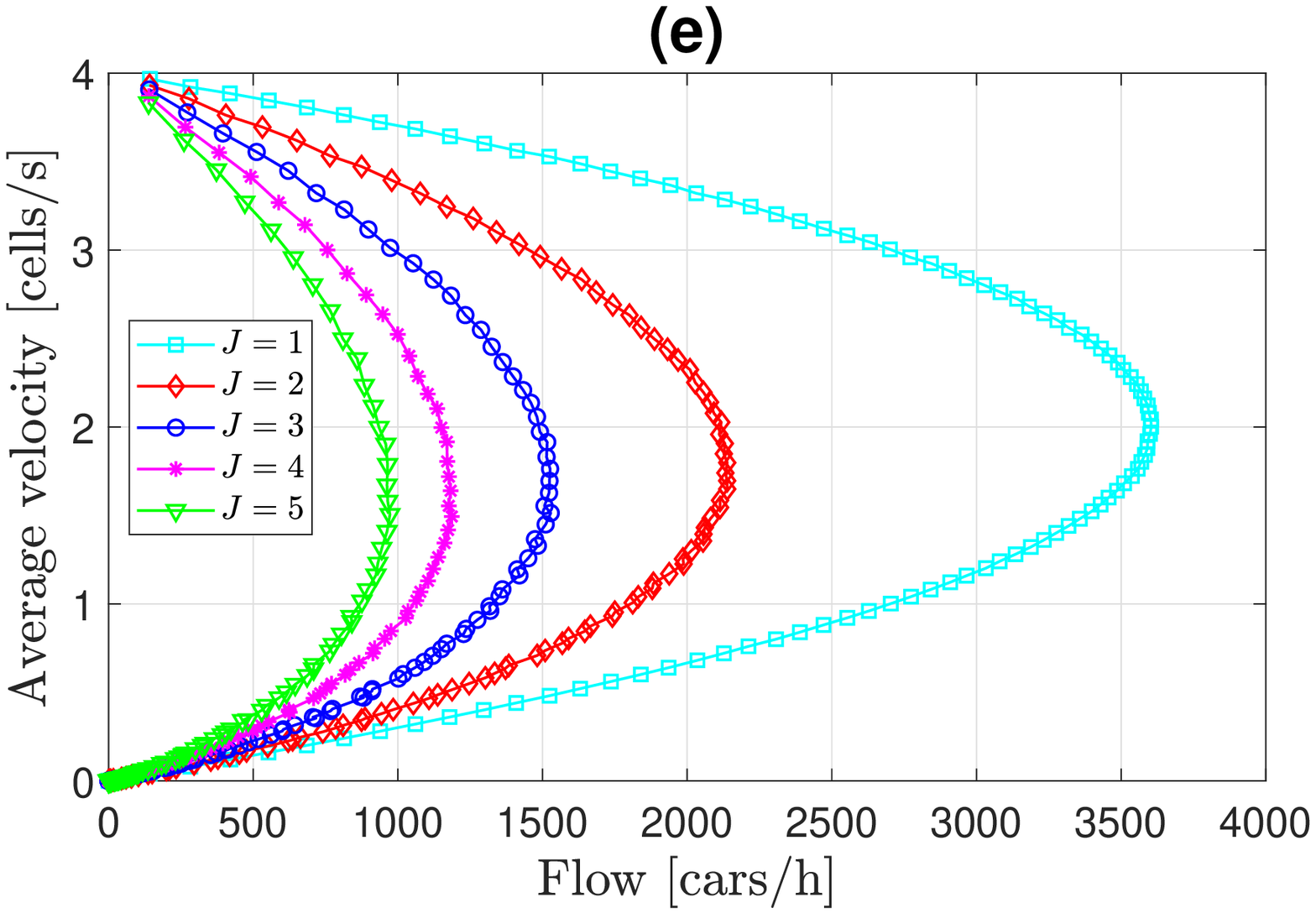} \hfill
\includegraphics[width=.47\textwidth]{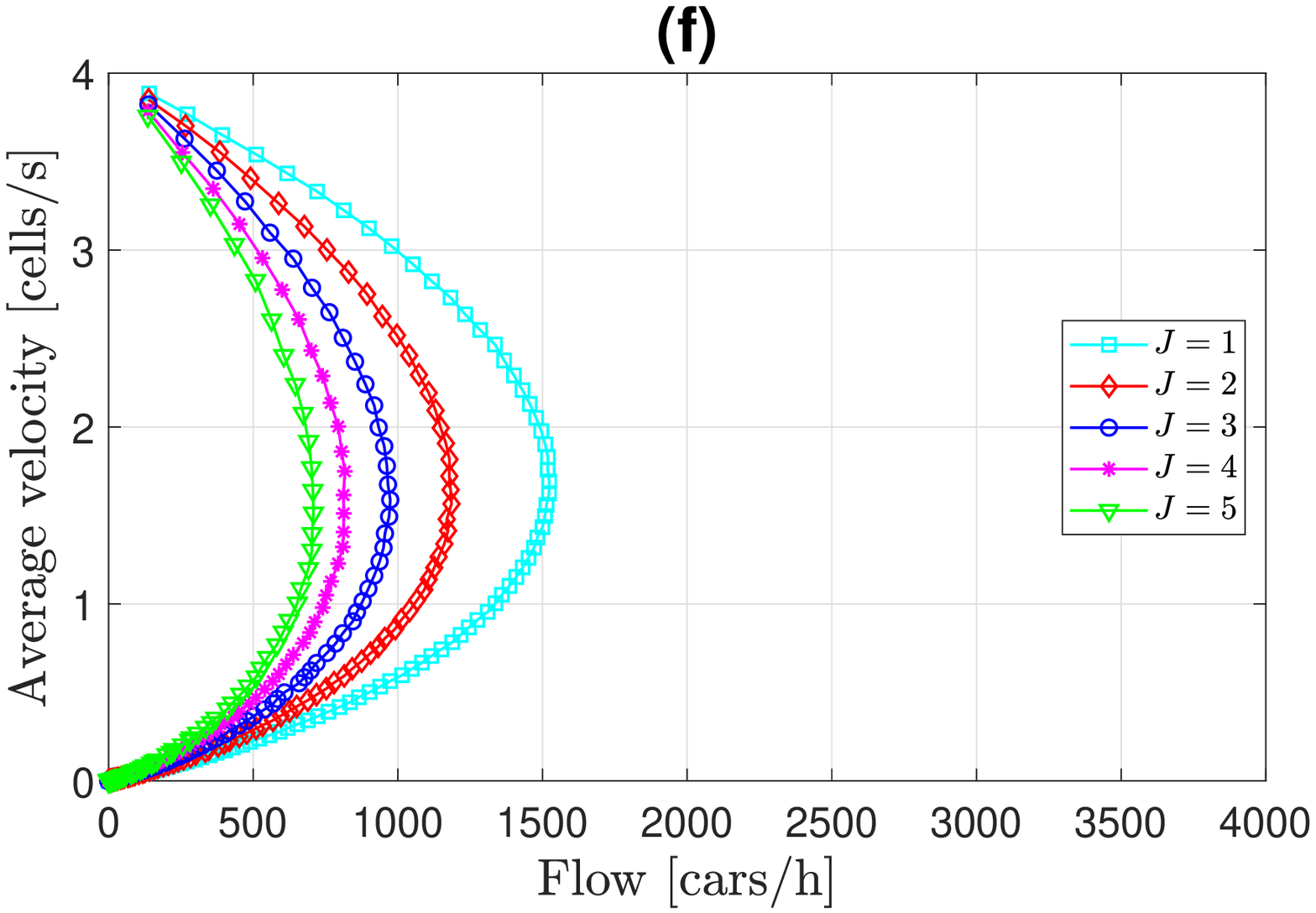}

\caption{Comparison results of the traffic flow on the one-lane highway with five different values of the multiple move parameter $\Jump$. (a)(b): Long-time averages of the density-flow relationship; (c)(d): Ensemble-averaged velocity of cars versus the density $\bar{\rho}$; (e)(f): Long-time averages of the flow-velocity relationship. (left panels): Results of $\Eo=10^4$. (right panels): Results of $\Eo=0.1$. Note that for each value of $\Jump=1,2,\ldots, 5$, the fluxes of the KMC simulations in (a) and (b) agree with the macroscopic averaged fluxes in \eqref{eq:fluxavginf} and \eqref{eq:genavg0} (shown as the dashed and dotted black curves, respectively). The slowdown  interaction function: $g(x) =(1-x)^2$ in \eqref{eq:g-2}. }
                                                           \label{fig:flowJ2}
\end{center}
\end{figure}

Fig.~\ref{fig:flowJ}(a) shows the density-flow relationship for the case of  $\Eo=10^4$ with $\Jump$ increasing from $1$ to $5$. The fluxes match beautifully with the macroscopic averaged fluxes \eqref{eq:fluxavginf} (shown as the dashed black curves). The case of $\Jump=1$ corresponds to the LWR type model, where the curve is symmetric and concave. We note that the same KMC results have also been shown for $\Eo=10^4$ and $\Jump=1$ in Fig.~\ref{fig:flowE}(a). For $\Jump\geq2$, the curves become neither convex nor concave, and have a right-skewed asymmetry.  Moreover, for a fixed $\bar{\rho}$, the magnitude of the flow $\Fbar$ decreases with increasing $\Jump$. The density-flow curves for $\Jump=1$ to $5$ take their maxima at the critical density $\rhoc=\frac12$ to $\frac16$, respectively, which is consistent with $\rho_c^\Eo\to\frac{1}{\Jump+1}$ as $\Eo\to\infty$ in \eqref{eq:rhocrit}. For the case of $\Eo=0.1$ shown in Fig.~\ref{fig:flowJ}(b), the microscopic fluxes agree with the macroscopic averaged fluxes \eqref{eq:fluxavg0} very well (shown as the dotted black curves). The critical density takes $\rhoc=\frac13$ to $\frac17$ for $\Jump=1$ to $5$, respectively, as $\rho_c^\Eo\to\frac{1}{\Jump+2}$ as $\Eo\to0$ in \eqref{eq:rhocrit}. Moreover, the averaged results of $\Eo=10^4$ for the cases $\Jump=2$ to $5$ in the left panels of Fig.~\ref{fig:flowJ} are statistically equal to the ones of $\Eo=0.1$ for the cases $\Jump=1$ to $4$ shown in the right panels of Fig.~\ref{fig:flowJ}.

The fundamental diagrams of the density-velocity relationship in Figs.~\ref{fig:flowJ}(c) and (d) show that at the same $\bar{\rho}$, the ensemble-averaged velocity $\vbar$ decreases as $\Jump$ increases. For the cases of $\Eo=10^4$ and $0.1$, the ensemble-averaged velocity of the KMC simulations agrees with limits \eqref{eq:fluxavginf} and \eqref{eq:fluxavg0}, (shown as the dashed and dotted black curves, respectively).

Figs.~\ref{fig:flowJ}(e) and (f) of the flow-velocity relationship show that for a fixed value of $\vbar$, the magnitude of the flow $\Fbar$ decreases with increasing $\Jump$. In Fig.~\ref{fig:flowJ}(e) of $\Eo=10^4$, as $\Jump$ increases, the maximum value of the flow decreases a lot from $\approx 3600$ cars per hour down to $\approx 965$ cars per hour, and the critical value $\vbar_{\rm crit}$ decreases a bit and becomes slightly lower than $2$ cells per second. For the case of $\Eo=0.1$ in Fig.~\ref{fig:flowJ}(f), the value of $\vbar_{\rm crit}$ does not change too much as $\Jump$ increases.

In Fig.~\ref{fig:flowJ2}, we show the three fundamental diagrams for the slowdown relation of $g(x)=(1-x)^2$  in \eqref{eq:g-2} with $\Eo=10^4$ and $\Eo=0.1$, respectively. Here, we observe the essentially same phenomena about the fundamental diagrams as shown in Fig.~\ref{fig:flowJ}. In particular, the averaged results for the case of $\Eo=10^4$ in the left panels of Fig.~\ref{fig:flowJ2} are statistically equal to the ones shown in the left panels of Fig.~\ref{fig:flowJ} as the long-time averaged flux and the averaged velocity still satisfy \eqref{eq:fluxavginf}. But for $\Eo=0.1$, the density-flow curves in Fig.~\ref{fig:flowJ2}(b) for all cases of $\Jump=1$ to $5$ take their maxima at the critical density $\rhoc=\frac14$ and $\frac18$, respectively, which is consistent with \eqref{eq:genavg0}.
Moreover, for all cases of $\Jump=1$ to $5$, the density-flow curves shown in Fig.~\ref{fig:flowJ2}(b) and the density-velocity curves in Fig.~\ref{fig:flowJ2}(d) are consistent with the macroscopic averaged flux and velocity in \eqref{eq:genavg0}, and these curves are below the corresponding curves shown in Fig.~\ref{fig:flowJ}(b) and (d), respectively. Meanwhile, the flow-velocity curves in Fig.~\ref{fig:flowJ2}(f) are left of the corresponding ones shown in Fig.~\ref{fig:flowJ}(f). As we pointed out in the previous section, this is because the function $g(x)=(1-x)^2$ introduces a stronger slow-down effect than $g(x)=(1-x)$ does.

\section{Conclusion}\label{sec:conclusion}

We have presented a class of one-dimensional cellular automata (CA) models on traffic flows, featuring nonlocal look-ahead interactions. We extended the Arrhenius type look-ahead rule in our previous work \cite{SuT20} to more general functions for characterizing the nonlocal slowdown effect. The look-ahead rule also features a novel idea of multiple moves, which plays a key role in recovering the right-skewed non-concave flux in the macroscopic dynamics. Through a semi-discrete mesoscopic stochastic process, we derived the coarse-grained macroscopic dynamics of the CA model.

To simulate the proposed CA models with general slowdown functions, we developed an accelerated KMC algorithm to improve computational efficiency. In the KMC method, the dynamics of cars is described in terms of the transition rates corresponding to possible configurational changes of the system, and then the corresponding time evolution of the system can be expressed in terms of these rates. For models with global look-ahead interactions, it is computationally costly to obtain all transition rates in each step due to its nonlocal nature. Therefore, in our accelerated KMC method, we take a new way to evaluate the transition rates by updating from the previous steps, which can reduce the
the cost from $\mathcal{O}(\NL^2)$ to $\mathcal{O}(\NL)$, where $\NL$ denotes the number of cells in the lattice.

The numerical experiments verified the computational efficiency of the accelerated KMC algorithm over the standard KMC method. Our numerical results show that the fluxes of the KMC simulations agree with the coarse-grained macroscopic averaged fluxes under various parameter settings. We obtained fundamental diagrams that display several important observed traffic states.

As one of our main goals is to compare the two look-ahead rules, we propose our CA models in a closed system and take the periodic boundary conditions to keep the number of cars and the density constant in a single simulation. Therefore, we have not applied our models to simulate some more complex non-stationary features, such as traffic breakdowns at bottlenecks \cite{Ker04}. It is possible to improve the models further in the following directions. We can include entrances and exits in the models by adding dynamical mechanisms such as adsorption/desorption. In reality, there are multi-lanes on highways and fast vehicles may change lanes to bypass slow ones. We also need to consider different types of vehicles, such as cars and trucks with unequal sizes and speeds. More complicated models addressing these aspects will be explored in the future.

\section*{Acknowledgment}
YS is partially supported by the NSF grants DMS-1913146 and DMS-1954532, a SC EPSCoR GEAR-CRP award (20-GC04), and a UofSC ASPIRE grant (2019-2021). CT is partially supported by the NSF grants DMS-1853001 and DMS-2108264, and a UofSC ASPIRE grant (2020-2022).


\end{document}